\documentclass[twoside, 11pt]{article}

\usepackage{amssymb}
\usepackage{amsfonts}
\usepackage{amsthm}
\usepackage{amsmath}
\usepackage{amscd}

\begin{document}
\pagestyle{myheadings}
\markboth{BASARAB}
{Right-angled Artin groups\/}
\title{{\huge On the arboreal structure of
right-angled Artin groups\/}}

\author{\c SERBAN A. BASARAB\\
Institute of Mathematics "Simion Stoilow" of the
Romanian Academy \\
P.O. Box 1--764\\
RO -- 70700 Bucharest 1, ROMANIA\\
\texttt{{\itshape e-mail}: Serban.Basarab@imar.ro}}
\date{} 

\maketitle
\vspace{8mm}

\newtheorem{te}{Theorem}[section]
\newtheorem{ex}[te]{Example}
\newtheorem{pr}[te]{Proposition}
\newtheorem{rem}[te]{Remark}
\newtheorem{rems}[te]{Remarks}
\newtheorem{co}[te]{Corollary}
\newtheorem{lem}[te]{Lemma}
\newtheorem{prob}[te]{Problem}
\newtheorem{exs}[te]{Examples}
\newtheorem{defs}[te]{Definitions}
\newtheorem{de}[te]{Definition}
\newtheorem{que}[te]{Question}
\newtheorem{aus}[te]{}


\newenvironment{aussage}%
\renewcommand{\theequation}{\alph{equation}}\setcounter{equation}{0}%

\def\R{\mathbb{R}}
\def\bbbm{{\rm I\!M}}
\def\N{\mathbb{N}}
\def\F{\mathbb{F}}
\def\H{\mathbb{H}}
\def\I{\mathbb{I}}
\def\K{\mathbb{K}}
\def\P{\mathbb{P}}
\def\D{\mathbb{D}}
\def\Q{\mathbb{Q}}
\def\Z{\mathbb{Z}}
\def\C{\mathbb{C}}
\def\T{\mathbb{T}}
\def\L{\mathbb{L}}
\def\cA{\mathcal{A}}
\def\cB{\mathcal{B}}
\def\cC{\mathcal{C}}
\def\cD{\mathcal{D}}
\def\cE{\mathcal{E}}
\def\cF{\mathcal{F}}
\def\cG{\mathcal{G}}
\def\cH{\mathcal{H}}
\def\cI{\mathcal{I}}
\def\cJ{\mathcal{J}}
\def\cK{\mathcal{K}}
\def\cL{\mathcal{L}}
\def\cM{\mathcal{M}}
\def\cN{\mathcal{N}}
\def\cO{\mathcal{O}}
\def\cP{\mathcal{P}}
\def\cQ{\mathcal{Q}}
\def\cR{\mathcal{R}}
\def\cS{\mathcal{S}}
\def\cT{\mathcal{T}}
\def\cU{\mathcal{U}}
\def\cV{\mathcal{V}}
\def\cW{\mathcal{W}}
\def\cX{\mathcal{X}}
\def\cY{\mathcal{Y}}
\def\cZ{\mathcal{Z}}
\def\epl{e^{+}}
\def\eb{e^{\bullet}}
\def\Ep{E^{+}}
\def\Eb{E^{\bullet}}


\newcommand{\notdiv}{{\not{|}\,}}
\newcommand{\ctg}{{\rm ctg\,}}
\newcommand{\sh}{{\rm sh\,}}
\newcommand{\ch}{{\rm ch\,}}
\newcommand{\Spek}{{\rm Spek\,}}
\newcommand{\Ker}{{\rm Ker\,}}
\newcommand{\cont}{\mbox{\rm cont\,}}
\newcommand{\rank}{\mbox{\rm rank\,}}
\newcommand{\codim}{\mbox{\rm codim\,}}
\newcommand{\ssqrt}[2]{\sqrt[\scriptstyle{#1}]{#2}}
\newcommand{\dbigcup}{\displaystyle\bigcup}
\newcommand{\dprod}{\displaystyle\prod}
\newcommand{\dbigcap}{\displaystyle\bigcap}
\newcommand{\doplus}{\displaystyle\bigoplus}
\newcommand{\dinf}{\displaystyle\inf}
\newcommand{\dint}{\displaystyle\int}
\newcommand{\dsum}{\displaystyle\sum}
\newcommand{\comb}[2]{\left(\begin{array}{c}#1\\#2\end{array}\right)}
\newcommand{\dlim}{\displaystyle\lim_{\stackrel{\longleftarrow}{}}}
\def\be{\begin{equation}}
\def\ee{\end{equation}}
\newcommand{\rs}{respectiv }
\newcommand{\G}{Z^1(\Gam,A)}

\def\bp{\begin{proof}}
\def\ep{\end{proof}}
\def\hf{\hfill $\square$}
\def\ben{\begin{enumerate}}
\def\een{\end{enumerate}}
\def\ba{\begin{eqnarray*}}
\def\ea{\end{eqnarray*}}

\def\ve{\varepsilon}
\def\ls{\leqslant}
\def\gs{\geqslant}
\def\lla{\longleftarrow}
\def\lra{\longrightarrow}
\def\Llra{\Longleftrightarrow}
\def\Lra{\Longrightarrow}
\def\p{\perp}
\def\wp{\,\widehat{1/p}}
\def\w4{\,\widehat{1/4}}
\newcommand{\wh}{\widehat}
\newcommand{\la}{\langle}
\newcommand{\ra}{\rangle}
\newcommand{\wt}{\widetilde}
\newcommand{\sm}{\setminus}
\newcommand{\sse}{\subseteq}
\newcommand{\es}{\varnothing} 

\newcommand{\f}{\frac}
\newcommand{\q}{\quad}
\newcommand{\n}{\vartriangleleft}


\newcommand{\al}{h }
\newcommand{\De}{\Delta}
\newcommand{\del}{\delta}
\newcommand{\eps}{\varepsilon}
\newcommand{\gam}{\gamma}
\newcommand{\Gam}{\Gamma}
\newcommand{\Lam}{\Lambda}
\newcommand{\lam}{\lambda}
\newcommand{\om}{\omega}
\newcommand{\Om}{\Omega}
\newcommand{\ovG}{\overline G}
\newcommand{\si}{\sigma}


\def\mA{\mathbb{A}}
\def\mC{\mathbb{C}}
\def\mN{\mathbb{N}}
\def\mQ{\mathbb{Q}}
\def\mR{\mathbb{R}}
\def\mZ{\mathbb{Z}}
\def\mF{\mathbb{F}}
\def\cA{\mathcal{A}}
\def\cC{\mathcal{C}}
\def\cD{\mathcal{D}}
\def\cE{\mathcal{E}}
\def\cF{\mathcal{F}}
\def\cG{\mathcal{G}}
\def\cM{\mathcal{M}}
\def\cP{\mathcal{P}}
\def\cS{\mathcal{S}}
\def\cO{\mathcal{O}}
\def\cL{\mathcal{L}}
\def\cQ{\mathcal{Q}}
\def\cK{\mathcal{K}}
\def\cH{\mathcal{H}}
\def\cX{\mathcal{X}}
\def\i{\rm i\/}
\begin{abstract}
The present article continues the study of median groups
initiated in \cite{Hyp, AC, Rav}. Some classes of median groups 
are introduced and investigated, with a stress upon the
class of the so called $A$-{\em groups} which contains
as remarkable subclasses the 
{\em lattice ordered groups} and the
{\em right-angled Artin groups}. Some general results concerning
$A$-groups are applied to a systematic study of the arboreal
structure of right-angled Artin groups. Structure theorems
for foldings, directions, quasidirections and centralizers are proved.

\smallskip


\smallskip

\end{abstract}

\bigskip

\section{Introduction}

$\quad$ In his paper \cite{AC}, an improved version of 
the preprints \cite{ArtCox1, ArtCox2}, the author applied 
the theory of generalized trees (median sets) developed in 
\cite{DF, Dual} to elaborate a general theory of median
(or arboreal) groups and use it to the investigation of 
a remarkable class of groups called {\em partially 
commutative Artin-Coxeter groups}.

Recall that by a \textit{median} or \textit{arboreal group} we mean a 
group $G$ endowed with a ternary operation $Y : G^3 \lra G$ making it a 
\textit{median set} or \textit{generalized tree} such that 
$uY(x,y,z) = Y(ux,uy,uz)$ for all $u,x,y,z \in G$. 

Equivalently, according to \cite[Proposition 2.2.1.]{AC}, 
by a median group we can also understand a group $G$ endowed 
with a meet-semilattice operation $\cap\,$, with associated 
order $\subset\,$, satisfying the following three axioms : 

\medskip

$(1)\,\forall x \in X, 1 \subset x$

$(2)\,\forall x, y , z \in X, 
x \subset y$ and $y \subset z \Lra z^{-1}y \subset z^{-1}x$

$(3)\,\forall x, y \in X, x^{-1}(x \cap y) \subset x^{-1}y$

In a median group $G$, an ordered pair $(x, y) \in G^2$ is said to
be {\em reduced} (write $xy = x \bullet y$) if $x \subset xy$, 
i.e. $x^{- 1} \cap y = 1$. For all $x, y \in G$, $x \cap y$ is
the unique element $z \in G$ satisfying $x = z \bullet (z^{- 1}x),
y = z \bullet (z^{- 1}y)$, and $x^{- 1}y = (x^{- 1}z) \bullet (z^{- 1}y)$.
Notice also that $xy \cap xz = x(y \cap z)$ (in particular,
$xy \subset xz \Llra y \subset z$) provided the pairs 
$(x, y)$ and $(x, z)$ are reduced. In other words, $x^{- 1}y \cap
x^{- 1}z = x^{- 1}(y \cap z)$ whenever $x \subset y$ and $x \subset z$.

The elements $x$ and $y$ of a median group $G$ are said to be 
\textit{orthogonal} (write $x \perp y$) if $x \cap y = 1$ and 
there exists the join $x \cup y$ (write $x \cup y \neq \infty$). 
According to \cite[2.2.]{AC}, a median group $G$ is said to be a 
$\perp$-\textit{group} if $\forall x, y \in G, x \perp y \Lra x \cup y = xy$, 
in particular, $xy = yx$. 

Notice that in a $\perp$-group $G$, 
for all $x, y \in G$, $x \cup y \neq \infty \Lra 
x \cup y = x(x \cap y)^{- 1}y = y(x \cap y)^{- 1}x$. 
In particular, a subgroup $H$ of $G$ is median, i.e.
$x, y, z \in H \Lra Y(x, y, z) \in H$, iff 
$x, y \in H \Lra x \cap y \in H$, while $H$ is convex,
i.e. $x, y \in H \Lra \forall z \in G, Y(x, y, z) \in H$, iff
$\forall x \in H, y \in G, y \subset x \Lra y \in H$.
By \cite[Corollary 2.8.]{Rav}, for any $\perp$-group $G$ there 
exists a canonical simple transitive action of $G$ on a 
subdirect product of locally linear median sets.

Given a group $G$ and a set $S \sse G$ of generators 
such that $1 \not\in S$ and 
$S_1 := S \cap S^{-1} = \{s \in S \mid s^2 = 1\,{\rm holds\,in}\,G\}$, 
it turns out by \cite[Theorem 2.4.1.]{AC} that the 
partial order on $G$ defined by $x \subset y \Llra
l(x) + l(x^{-1}y) = l(y)$, where $l : G \lra \N$ denotes the
canonical length function on $(G, S)$, makes $G$ a $\perp$-group if and 
only if $(G, S)$ is a {\em partially commutative Artin-Coxeter system},
i.e. the group $G$ admits the presentation 
$$G = <S; s^2 = 1\, {\rm for}\, s \in S_1, [s, t] = 1 
\,{\rm for}\, s,t \in S, s \neq t,\,
st = ts\,{\rm holds\, in}\, G>$$
Thus the partially commutative Artin-Coxeter groups are 
identified with the simplicial $\perp$-groups.

The present paper, an improved version of the preprint \cite{ArtCox3}, 
is devoted to a systematic study of the arboreal 
structure of the systems $(G,S)$ above which satisfy the additional
restrictive assumption that $S \cap S^{-1} = \emptyset$. 
Such systems were introduced by Baudisch in \cite{Bau1, Bau2} under
the name of \textit{semi-free groups}, and extensively studied
in the last years under various names (\textit{right-angled
Artin groups, free partially commutative groups, graph groups})
by people working in combinatorial and geometric group theory, 
associative algebras, computer science (see for instance the
long bibliography to the survey article \cite{Charney}). 
Their nice properties were exploited by 
Bestvina and Brady \cite{BB} in their construction of examples of 
groups which are of type (\textit{FP}) but are not finitely
presented, as well as by Crisp and Paris \cite{Crisp} in their proof of
a conjecture of Tits on the subgroup generated by the squares of the
generators of an Artin group. According to \cite{Davis}, the finitely 
generated right-angled Artin groups (moreover, the weakly
partially commutative Artin-Coxeter groups as defined in 
\cite[1.1.]{AC}) are linear and hence equationally noetherian.

The outline of the paper is as follows. Some notions and 
basic facts from \cite{DF} on \textit{congruences} and \textit{quasidirections} 
on median sets are recalled in Section 2. Some classes of median groups
are introduced and investigated in Section 3. Amongst them, the class
of the so called $A$-{\em groups} contains as remarkable subclasses the
$l$-groups, not necessarily commutative, and the right-angled Artin groups.
The {\em cyclically reduced elements} of $A$-groups are studied 
in Section 4, while Sections 5 and 6 are devoted to the 
main properties of the preorders 
$\mathop{\preceq}\limits_{w}$ and the foldings $\varphi_w$ 
naturally associated to any element $w$ of an $A$-group. 

The general theory of $A$-groups is further applied in the last two sections 
of the paper to the particular case of right-angled Artin groups. One shows 
in Section 7 that in a right-angled Artin group $G$ the preorders 
$\mathop{\preceq}\limits_{w}$ determine quasidirections 
$\mathop{\bullet}\limits_{w}$ which are described as limits 
of sequences of operators $w^n\varphi_w$ for 
$n \lra \infty$.

The main results of the paper contained in Section 8 provide 
descriptions of the quasidirections $\mathop{\bullet}\limits_{w}$, 
the foldings $\varphi_w$ and the centralizers $Z_G(w)$ in terms 
of the corresponding invariants $\mathop{\bullet}\limits_{p}$, 
$\varphi_p$ and $Z_G(p)$, where $p$ ranges over a finite set 
$Prim(w)$ of \textit{primitive elements} canonically associated 
to any element $w$ of a right-angled Artin group $G$.

\bigskip


\section{Congruences and quasidirections on median sets}

$\quad$ In this section having a preliminary character we recall 
some notions and basic facts from \cite{DF} which will be used later.

By a {\it median set} or {\it generalized tree} we understand 
a set $X$ endowed with a ternary operation $Y : X^3 \lra X$, 
called {\it median}, satisfying the following equational axioms:

(i) {\it Symmetry} : $Y(x, y, z) = Y(y, x, z) = Y(x, z, y)$,

(ii) {\it Absorptive law} : $Y(x, y, x) = x$, and

(iii) {\it Selfdistributive law} : $Y(Y(x, y, z), u, v) =
Y(Y(x, u, v), y, Y(z, u, v))$.

\medskip

In a median set $X$, for any $a \in X$, the binary operation 
$(x, y) \mapsto x \mathop {\vee }\limits_{a} y := Y(x, a, y)$
makes $X$ a join-semilattice with the last element $a$; let
$\mathop {\leq}\limits_{a}$ denote the corresponding partial
order.

A subset $I$ of a median set $X$ is said to be {\it convex} if
$Y(x, y, z) \in I$ whenever $x, y \in I, z \in X$. As the
intersection of an arbitrary family of convex subsets is also convex,
we may speak on the {\it convex closure} of any subset $A$ of
$X$ and denote it by $[A]$. In particular, set 
$[a, b] =: [\{ a, b \} ]$ for $a, b \in X$.

By a {\it cell} of a median set $X$ we mean a convex subset $C$
of $X$ for which there are $a, b \in X$ such that $C = [a, b]$. 
Given a cell $C$, every element $a \in X$ for which there exists $b \in X$
such that $C = [a, b]$ is called an {\it end} of $C$. The
(non-empty) subset of all ends of a cell $C$, denoted by 
$\partial C$ and called the {\it boundary} of $C$, is a 
median subset of $C$, and the mapping $\neg$ assigning to each
$a \in \partial C$ the unique end $\neg \ a$ of $C$ for which
$C = [a, \neg\, a]$ is an involutory automorphism of the median
set $\partial C$. Note also that for a given $a \in \partial C$,
the cell $C$ becomes a bounded distributive lattice with respect to  
the order $\mathop {\leq}\limits_{a}$, with the join 
$\mathop {\vee}\limits_{a}$, the meet 
$\mathop {\vee}\limits_{\neg \ a}$, the last element $a$, 
and the least element $\neg \, a$, while its boundary $\partial C$
is identified with the boolean subalgebra consisting of 
those elements which have (unique) complements.

The median set $X$ is called {\it locally boolean}, resp.
{\it locally linear}, if $C = \partial C$ for every cell $C$
of $X$, resp. $\partial [x, y] = \{ x, y \}$ for all $x, y \in X$.
$X$ is called {\it simplicial} (or {\it discrete} or 
{\em locally finite}) if every cell
of $X$ has finitely many elements. A  graph-theoretic 
description for simplicial median sets is given in 
\cite[Lemma 7.1., Proposition 7.3.]{UC1} In particular, 
the {\it trees}, i.e. the acyclic connected graphs, are 
naturally identified with the simplicial locally linear median sets.

Note that the convex closure of a finite subset of a simplicial
median set is finite too, and hence every non-empty convex subset is 
retractible. To any simplicial median set $X$ one assigns 
an integer-valued "distance" function 
$d : X \times X \to \N$, where for
$x, y \in X, d(x, y)$ is the length of some (of any) maximal
chain in the finite distributive lattice 
$([x, y], \mathop{\leq}\limits_{y})$. With respect to $d$, 
$X$ becomes a $\Z$-metric space such that for all 
$x, y \in X, [x, y] = \{z \in X \mid d(x, z) + d(z, y) = d(x, y)\}$, 
and the mapping $[x, y] \to [0, d(x, y)], z \mapsto d(x, z)$, 
induced by $d$, is onto. In particular, $X$ is a tree iff 
for all $x, y \in X$, the mapping above is bijective.

\medskip

\subsection{Congruences on median sets}

$\quad$ Given a median set $X$, a {\it congruence} on $X$ is 
an equivalence relation $\rho$ on $X$ which is compatible 
with the median $Y$, i.e. for all $a, b, x, y \in X$, 
$$(x, y) \in \rho \Lra (Y(a, b, x), Y(a, b, y)) \in \rho.$$

The congruences on $X$ form a lattice $Cong (X)$ with 
a least and a last element under the inclusion of relations. 
Moreover, according to \cite[Proposition 1.6.1.]{DF}, the lattice 
$Cong (X)$ is a Heyting algebra (in particular, 
a bounded distributive lattice), i.e. for every pair $(\rho, \gam)$ 
of congruences on $X$ there exists a unique congruence 
$\mu:=\rho \to \gam$ subject to $\theta \subseteq \mu \Llra
\theta \cap \rho \sse \gam$ 
for all $\theta \in Cong(X)$, namely the congruence
$$\mu = \{(a,b) \in X \times X \mid \forall x, y 
\in [a,b], (x,y) \in \rho \Lra (x,y) \in \gamma\}.$$

In particular, for $\gamma = \Delta$, the equality on $X$, 
we obtain the {\it negation} of $\rho$
$$\neg \rho : = \rho \to \Delta = \{(a,b) \in X \times X \mid 
\rho \mid_{[a,b]} = \Delta|_{[a,b]}\}.$$ 

By \cite[Corollary 1.6.2.]{DF}, $Cong(X)$ is a boolean algebra 
provided the median set $X$ is simplicial. 

Given a simplicial median set $X$ and a congruence $\sim$ 
on $X$, let $\equiv$ denote the negation (the complement) 
$\neg \sim$ of $\sim$ in the boolean algebra $Cong (X)$. 
For every $a \in X$, set $\widetilde a = \{x \in X \mid x \sim a\}$, 
$\mathop{a}\limits^{\equiv} = \{x \in X \mid x \equiv a \}$, 
and let $\varphi_a$, resp. $\psi_a$, denote the 
folding induced by the (retractible) convex subset 
$\widetilde a$, resp. $\mathop{a}\limits^{\equiv}$. Thus for all 
$x \in X, [a, x] \cap \widetilde a = [a, \varphi_a(x)]$ and 
$[a,x] \cap \mathop{a}\limits^{\equiv} = [a, \psi_a(x)]$.

\begin{lem} Let $\sim$ be a congruence on a simplicial median set $X$, 
with its negation $\equiv$. Then, the following assertions are equivalent.

$(1)\,$ For all $a, b \in X$, the intersection 
$\widetilde a \cap \mathop{b}\limits^{\equiv}$ is nonempty.

$(2)\,$ For all $a, b \in X, \varphi_a(b) \equiv b$, 
i.e. for all $a \in X$, the embedding $\widetilde a \lra X$ 
induces a median set isomorphism $\widetilde a \lra X/\equiv$.

$(3)\,$ For all $a, b \in X, \psi_a(b) \sim b$, i.e. for all
$a \in X$, the embedding $\mathop{a}\limits^{\equiv} \lra X$ 
induces a median set isomorphism 
$\mathop{a}\limits^{\equiv} \lra X/\sim$.

$(4)\,$ For all $a, b \in X, \varphi_a(b) = \psi_b(a)$.

$(5)\,$ For all $a, b \in X, [a, b] = [\varphi_a(b), \psi_a(b)]$.

$(6)\,$ For every quasi-linear cell $[a, b]$ (i.e. 
$\partial [a, b] = \{a, b\}$), either $a \sim b$ 
or $a \equiv b$.

$(7)\,$ For every cell $[a, b]$ with three elements, 
either $a \sim b$ or $a \equiv b$.
\end{lem}

\bp 
The implications $(4) \Lra (2), (4) \Lra (3), (2) \Lra (1), 
(3) \Lra (1)$ and $(6) \Lra (7)$ are trivial.

$(1) \Lra (4)$. If $\widetilde a \cap 
\mathop{b}\limits^{\equiv}$ is non-empty, then obviously 
$\widetilde a \cap \mathop{b}\limits^{\equiv} = \{c\}$ 
is a singleton. Note also that $c = Y(c, c, b) \sim Y(a, c, b)$ 
and $c = Y(a, c, c) \equiv Y(a, c, b)$, therefore $c = Y(a, c, b)$, 
i.e. $c \in [a, b]$. Consequently, 
$\{c\} = [a, b] \cap \widetilde a \cap \mathop{b}\limits^{\equiv} 
= [a, \varphi_a(b)] \cap [b, \psi_b(a)] = [Y(a, \varphi_a(b),b) = 
\varphi_a(b), Y(a, \varphi_a(b), \psi_b(a)] = [Y(a, b, \psi_b(a)) = 
\psi_b(a), Y(\varphi_a(b), b, \psi_b(a))]$, and hence 
$c = \varphi_a(b) = \psi_b(a)$ as required.

$(4) \Lra (5)$. The inclusion 
$[\varphi_a(b), \psi_a(b)] \subseteq [a,b]$ is obvious. 
On the other hand, $Y(a,\varphi_a(b), 
\psi_a(b)) \sim Y(a,a, \psi_a(b)) = a$, and 
$Y(a, \varphi_a(b), \psi_a(b)) \equiv Y(a, \varphi_a(b), a) = a$, 
therefore $Y(a, \varphi_a(b), \psi_a(b)) = a$, i.e. 
$a \in [\varphi_a(b), \psi_a(b)]$. By symmetry, we get 
$b \in [\varphi_b(a), \psi_b(a)] = [\psi_a(b), \varphi_a(b)]$ 
(by assumption). Thus $[a, b] \sse [\varphi_a(b), \psi_a(b)]$ 
as desired.

$(5) \Lra (6)$. Since the cell $[a, b] = [\varphi_a(b), 
\psi_a(b)]$ is assumed to be quasilinear, it follows that either 
$\varphi_a(b) = b$, i.e. $a \sim b$, or $\psi_a(b) = b$, 
i.e. $a \equiv b$.

$(7) \Lra (2)$. As $\varphi_a = \varphi_{\varphi_a(b)}$, 
we may assume without loss that $\varphi_a(b) = a$, 
so we have to show that $a \equiv b$. We argue by 
induction on the "distance" $d := d(a, b)$. Since 
the cases $d = 0$ and $d = 1$ are trivial, we may assume 
$d \geq 2$. Let $c \in [a, b]$ be such that $d(c, b) = 2$, 
and let $e \in [c, b] \backslash \{c, b\}$. Since 
$\varphi_a(e) = \varphi_a(Y(a,e,b)) = Y(\varphi_a(a),e, 
\varphi_a(b)) = Y(a, e, a) = a$ and  $d(a, e) = d - 1 < d$, 
it follows by the induction hypothesis that $a \equiv e$, 
therefore $c \equiv a \equiv e$, as $c \in [a,e]$. We 
distinguish the following two cases:

Case $(i)\,$ : The cell $[c, b]$ has three elements, 
i.e. $[c, b] = \{c, e, b\}$. By assumption either 
$c \equiv b$ or $c \sim b$. In  the former case 
we get $a \equiv b$, as required, while in the 
latter case it follows that $c \sim e$ as $e \in [c, b]$. 
Since, on the other hand, $c \equiv e$, we get $c = e$, 
i.e. a contradiction.

Case $(ii)\,$ : The cell $[c, b]$ has four elements, say 
$[c, b] = [e, f] = \{c, b, e, f\}$. As we already know 
that $e \equiv a \equiv f$, we get $a \equiv b$ since 
$b \in [e, f]$.
\ep 

\begin{co} Given two complementary 
congruences $\sim$ and $\equiv$ on a simplicial median 
set $X$, assume that for each cell $[a,b]$ with three 
elements either $a \sim b$ or $a \equiv b$. For 
all $a \in X$, let $\varphi_a$, resp, $\psi_a$, 
denote the folding of $X$ induced by the convex subset 
$\widetilde a$, resp $\mathop{a}\limits^{\equiv}$. 
Then, for all $a \in X$, the median set morphism 
$X \lra \widetilde a \times \mathop{a}\limits^{\equiv}, x \mapsto 
(\varphi_a(x), \psi_a(x))$ is an isomorphism, whose 
inverse sends a pair 
$(x, y) \in \widetilde a \times \mathop{a}\limits^{\equiv}$ to
$\psi_x(y) = \varphi_y(x)$.
\end{co}

\medskip

\subsection{Directions and quasidirections on median sets}

$\quad$ By a {\it quasidirection} on a median set $X$ we 
understand a binary operation $\bullet$ on $X$ satisfying
the following four conditions :

i) $(X, \bullet)$ is a {\it band}, i.e. a semigroup 
in which all elements are idempotent,

ii) $a \bullet b \bullet c = a \bullet c \bullet b$ 
for all $a, b, c \in X$,

iii) for all $a \in X$, the left translation 
$X \to X, x \mapsto a \bullet x$ is a folding, i.e. \\
$a \bullet Y(x, y, z) = Y(a \bullet x, y, a \bullet z)$ 
for all $x, y, z \in X$, and

iv) for all $x, y, z \in X, Y(x, y, z) \bullet x = 
Y(x, y, z \bullet x)$.

Moreover, by \cite[Lemma 3.3.]{DF}, a stronger form of iv),
the symmetrical version of iii), is also satisfied :

iii)' for all $a \in X$, the right translation 
$X \to X, x \mapsto x \bullet a$ is a folding of $X$.

A quasidirection $\bullet$ on $X$ is said to be a 
{\it direction} if the band $(X, \bullet)$ is a 
semilattice, i.e. $x \bullet y = y \bullet x$ 
for all $x, y \in X$. In this case, iv) becomes 
supperflous. Any element $a$ of $X$ determines a 
direction $\mathop{\vee}\limits_{a}$ on $X$ given 
by $x \mathop{\vee}\limits_{a} y: = Y(x, a, y)$. 
Such directions are called {\it internal} or 
{\it closed}, while the other ones, if exist, are 
called {\it external} or {\it open}.

Call {\it trivial} the quasidirection defined by 
the rule $x \bullet y = x$.

A median set $X$ endowed with a quasidirection, 
resp. a direction, is said to be {\it quasidirected}, 
resp. {\it directed}.

According to \cite[Proposition 3.7.]{DF}, the mapping 
assigning to a binary operation $\bullet$ on the 
median set $X$ the binary relation 
$a \mathop{\leq}\limits_{l} b \Llra b \bullet a = b$ 
maps bijectively the set of quasidirections on $X$ 
onto the set of the preorders $\preceq$ on $X$ satisfying

i) $\preceq$ is compatible with the median of $X$, 
i.e. $\forall a, b, x, y \in X, x \preceq y \Lra
Y(a, b, x) \preceq Y(a, b, y)$; let $\sim$
denote the congruence induced by the preorder $\preceq$, 
and let $\equiv$ be its negation in the Heyting algebra $Cong\, (X)$;

ii) for all $a, b \in X$ there exists $c \in X$ such that 
$a \preceq c, b \preceq c$, and $a \equiv c$.

The inverse of the bijection above sends a preorder 
$\preceq$ as above to the quasidirection $\bullet$ given by 
$a \bullet b = Y(a, b, c)$ for some (for all) $c \in X$ 
subject to ii). Note also that $a \equiv b \Llra
a \bullet b = b \bullet a$.

The bijection above induces by restriction a bijection 
of the set of directions on $X$ onto the set of the 
orders of $X$ which are compatible with the median of 
$X$ such that any pair $(a, b)$ of elements in $X$ 
is bounded above.

According to \cite[Corollary 3.5.]{DF}, the canonical embedding 
$X \to X/\sim \times X/\equiv$ yields a representation 
of the quasidirected median set $(X, \bullet)$ as a 
subdirect product of a pair consisting of a directed 
median set $X/\sim$ and a trivially quasidirected median 
set $X/\equiv$, in such a way that the product 
$X/\sim \times X/\equiv$ is the convex closure 
of its median subset $X$.

Given a median set $X$, let $Dir\, (X), Fold\, (X)$ and 
$Q dir\, (X)$ respectively denote the set of directions, 
of foldings and of quasidirections on $X$. By \cite[Sections 8, 10]{DF},
$Q dir\, (X)$ becomes a directed median set with the median 
$(q_1, q_2, q_3) \mapsto Y(q_1, q_2, q_3)$ given by 
$a \mathop{\bullet}\limits_{Y(q_1, q_2, q_3)} b = Y(a 
\mathop{\bullet}\limits_{q_1} b, a \mathop{\bullet}\limits_{q_1} b, 
a \mathop{\bullet}\limits_{q_3} b)$, and the direction 
induced by the order $q_1 \leq q_2$ iff the preorder 
$\mathop{\preceq}\limits_{q_2}$ associated to $q_2$ is 
finer than $\mathop{\preceq}\limits_{q_1}$. 
The subset $Dir\, (X)$ is a median subset of $Q dir\, (X)$ 
consisting of the minimal elements under the order 
$\leq$ on $Q dir\, (X)$, while the injective mapping 
$X \to Dir\, (X), a \mapsto \mathop{\vee}\limits_{a}$, 
identifies $X$ with a convex subset of $Dir\, (X)$. 
On the other hand, by \cite[Proposition 9.1.]{DF}, $Fold\, (X)$ 
has a canonical structure of directed median set with 
the median defined by $Y(\varphi_1, \varphi_2, 
\varphi_3)(x) = Y(\varphi_1(x), \varphi_2(x), \varphi_3(x))$, 
and the direction induced by the order $\varphi \leq \Psi$ 
iff $\varphi(X) \subseteq \Psi(X)$. Note that the injective 
mapping $X \to Fold\, (X), a \mapsto (x \mapsto a)$ 
identifies $X$ with a median subset of $Fold\, (X)$.

According to \cite[Theorem 9.3.]{DF}, the mapping 
$\alpha : Dir\, (Fold\, (X)) \lra Fold\, (Dir\, (X))$, 
given by $a \mathop{\vee}\limits_{\alpha(d)(D)} b = 
(a \mathop{\vee}\limits_{d} b)(a \mathop{\vee}\limits_{D} b) =$ 
the value in a $\mathop{\vee}\limits_{D} b$ 
of the folding $a \mathop{\vee}\limits_{d} b$, for
$d \in Dir\, (Fold\, (X)), D \in Dir\, (X), a, b \in X$, 
is an isomorphism of median sets, 
while by \cite[Theorem 10.1.]{DF}, the map 
$\gamma : Fold\, (Dir\, (X)) \lra Q\, dir\, (X)$, 
given by $a \mathop{\bullet}\limits_{\gamma(\varphi)} b = 
a \mathop{\vee}\limits_{\varphi(a)} b$ for 
$\varphi \in Fold\, (Dir\, (X)), a, b \in X$, is an 
isomorphism of directed median sets.

Given a simplicial median set $X$ and a preorder 
$\preceq$ on $X$ which is compatible with the median $Y$, 
let $\sim$ denote the congruence induced by 
$\preceq$, with its complement $\equiv$ in the boolean 
algebra $Cong\, (X)$. Recall that $x \equiv y \Llra
\forall u, v \in [x, y], u\sim v \Lra u = v$. Assume 
that any pair $(a, b)$ of elements in $X$ is bounded above 
with respect to the preorder $\preceq$, i.e. there exists 
$c \in X$ such that $a \preceq c$ and $b \preceq c$. For 
$a, b \in X$, set $U_{a, b} = \{x \in [a, b] \mid a \preceq x \; 
{\mbox{and}} \; b \preceq x\}$. By assumption, the finite set 
$U_{a, b}$ is nonempty. Indeed, if $c$ is a common upper 
bound of the elements $a$ and $b$, then 
$a = Y(a, b, a) \preceq Y(a, b, c)$, and 
$b = Y(a, b, b) \preceq Y(a, b, c)$, therefore 
$Y(a, b, c) \in U_{a, b}$. Define the binary 
operation $\bullet$ on $X$ by 
$a \bullet b = \mathop{\vee}\limits_{a} U_{a,b}$. 

With the notation and the data above we have

\begin{lem}

$(1)\,$ $U_{a,b} = [a \bullet b, b \bullet a]$.

$(2)\,$ $a \bullet b \sim b \bullet a$,

$(3)\,$ $a \mathop{\leq}\limits_{b} a \bullet b \mathop{\leq}\limits_{b} b 
\bullet a \mathop{\leq}\limits_{b}  b$.

$(4)\,$ $a \preceq b \Llra b \bullet a = b$.
\end{lem}

\bp 
$(1)\,$ Since by assumption the preorder $\preceq$ is 
compatible with the median $Y$, it follows that the 
nonempty set $\{x \in X \mid a \preceq x$ and 
$ b \preceq x\}$ is a convex subset of $X$. 
Consequently, its image $U_{a, b}$ through the folding 
$X \longrightarrow X, x \mapsto Y(a, b, x)$, is a convex 
subset of the cell $[a, b]$. In particular, the cell 
$[a \bullet b, b \bullet a]$ is contained in $U_{a, b}$. 
On the other hand, for any $c \in U_{a, b}$ it follows by 
definition that $a \bullet b \in [a, c]$ and 
$b \bullet a \in [b, c]$. Consequently, 
$c \in [a \bullet b, b \bullet a]$ since otherwise,
by \cite[Corollary 5.2.2.]{Dual}, there exists 
a prime convex subset $P$ of $X$ such that 
$a \bullet b \in P,\, b \bullet a \in P$, and $c \not \in P$, 
therefore $c \in [a, b] \sse P$, a contradiction.

The statements $(2), (3)\,$ and $(4)$ are obvious.
\ep

The next lemma provides a characterization of those preorders 
on a simplicial median set which induce quasidirections.

\begin{lem} Let $\preceq$ be a preorder on 
a simplicial median set $X$ which is compatible with the 
median $Y$, such that any pair $(a, b)$ of elements of $X$ 
is bounded above with respect to $\preceq$. With the 
notation above, the following assertions are equivalent.

$(1)\,$ The binary operation $\bullet$ induced by the preorder 
$\preceq$ is a quasidirection on $X$.

$(2)\,$ $\forall a, b \in X, a \bullet b = b \Lra a \equiv b$.

$(3)\,$ For all $a, b, c \in X$ such that $[a, b] = \{a, c, b\}$ 
and $c \not\in \{a, b\}, c \preceq b \Lra a \preceq b$.
\end{lem}

\bp
$(1) \Lra (3)$. Let $a, b, c \in X$ be such 
that $[a, b] = \{a, c, b\}, c \not\in \{a, b\}$ and 
$c \preceq b$. Assuming that $a \not\preceq b$ it follows that 
$U_{a, b} = \{a\}$, i.e. $a \bullet b = b \bullet a = a$. 
In particular, $b \sim c$ since $c \preceq b, c \in [b, a]$ 
and $b \preceq a$. Since by assumption the binary operation 
$\bullet$ is a quasidirection on $X$, we get 
$a \equiv b$, contrary to $b \sim c, b \neq c$.

$(3) \Lra (2)$. Let $a, b \in X$ be such that 
$a \bullet b = b$, i.e. $U_{a, b} = \{b\}$. To show that 
$a \equiv b$ we argue by induction on the distance 
$d := d(a, b)$. Since the cases $d = 0$ and $d = 1$ are 
trivial, we may assume $d \geq 2$. Let $c \in [a, b]$ 
be such that $d(a, c) = 2$, and let 
$x \in [a, c] \backslash \{a, c\}$. It follows that 
$U_{x, b} = \{b\}$, i.e. $x \bullet b = b$, since 
$y \in U_{x, b} \Lra a \preceq b \preceq y$, 
and hence $y \in U_{a, b} = \{b\}$. As $d(x, b) = d - 1 < d$, 
it follows by the induction hypothesis that $x \equiv b$. 
We distinguish the following two cases :

Case $(i)\,$. The cell $[a, c]$ has four elements, say 
$[a, c] = [x, y] = \{a, c, x, y\}$. As we already know that 
$x \equiv b$ and $y \equiv b$, it follows that $a \equiv b$ 
since $a \in [x, y]$ and  $\equiv$ is a congruence on the 
median set $X$.

Case $(ii)\,$ The cell $[a, c]$ has three elements, say 
$[a, c] = \{a, x, c\}$. Assuming $a \sim x$ it follows by 
the assumption $(3)\,$ that $c \preceq a$, therefore $c \preceq x$. 
On the other hand, $x \preceq c$ since $c\in [x, b]$ and 
$x \preceq b$. Thus $x \sim c$, contrary to $x \equiv b, 
c \in [x,b], c \neq x$. Consequently, $a \not\sim x$, 
therefore $a \equiv x$ since $d(a,x) = 1$, and hence 
$a \equiv b$ as required.

$(2) \Lra (1)$. To conclude that the binary 
operation $\bullet$ induced by the preorder $\preceq$ 
is a quasidirection on $X$, by \cite[Proposition 3.7.]{DF} it 
suffices to show that $a \bullet b \equiv a $ for all 
$a, b \in X$. Thanks to the assumption $(2)\,$ we have to 
check the identity $a \bullet (a \bullet b) = a \bullet b$, 
i.e. $U_{a, a \bullet b} = \{a \bullet b\}$. Obviously, 
$a \bullet b \in U_{a, a \bullet b}$ since 
$a \preceq a \bullet b$. On the other hand, for any 
$c \in  U_{a, a \bullet b}$ we get $a \preceq c$ and 
$b \preceq a \bullet b \preceq c$, and hence 
$c \in U_{a, b}$. 
Since $c \in  U_{a, a \bullet b} 
\subseteq [a, a \bullet b]$, 
and $a \bullet b = \mathop{\vee}\limits_{a} U_{a, b} \in [a, c]$, 
it follows that $c = a \bullet b$ as required. \ep

\begin{rem} \em 
Consider the simplicial tree $X$ with three vertices 
$a, b, c$ and two geometric edges $(a, c)$ and $(c, b)$.
Let $\preceq$ be the complement in $X \times X$
of the subset $\{(a,b), (a,c)\}$. The relation $\preceq$ 
is a preorder on $X$ which is compatible with the tree 
structure, and any pair of elements of $X$ is 
bounded above with respect to $\preceq$. However the 
preorder $\preceq$ does not induce a quasidirection 
on $X$ since the condition $(3)\,$ above is not satisfied. 
Indeed, $c \preceq b$ but $a \not\preceq b$. Notice
that in this simple case, $Q dir(X) \cong Fold(X)$ is 
naturally identified with the directed median set 
$$\{\{a\}, \{b\}, \{c\}, [a, c] = \{a, c\}, [b, c] =
\{b, c\}, X = [a, b] = \{a, b, c\}\}$$ 
of cardinality $6$, consisting of the nonempty
convex subsets of $X$.
\end{rem}

\bigskip

\section{Some classes of median groups}

$\quad$ As shown in \cite{AC}, the class of partially
commutative Artin-Coxeter groups is naturally embedded into
a larger class of median groups consisting of the so
called $\perp$-groups, as defined in Introduction. Since
the right-angled Artin groups form a proper subclass of
the partially commutative Artin-Coxeter groups, it is 
natural to look for a proper subclass of the $\perp$-groups
which is adequate for the investigation of the arboreal
structure of right-angled Artin groups.

First of all notice that any {\it l}-group 
$(G, \mathop{.}, \leq, \wedge, \vee)$, not necessarily commutative, 
has a canonical structure of median group. Indeed, as the underlying 
lattice of $G$ is distributive, $G$ has a canonical structure of 
median set with the median defined by 
$Y(x, y, z) = (x \wedge y) \vee (y \wedge z) \vee (z \wedge x) =
(x \vee y) \wedge (y \vee z) \wedge (z \vee x)$. Obviously, the
median operation is compatible with the multiplication, so
$G$ becomes a median group. Notice that $x \subset y$ 
iff $x_+ \leq y_+$ and $x_- \leq y_-$, $(x \cap y)_+ = x_+ \wedge y_+,
(x \cap y)_- = x_- \wedge y_-$, where $x_+ = x \vee 1, x_- =
(x^{-1})_+ = (x \wedge 1)^{-1}$. 

For $x, y \in G, x \perp y$ iff $x$ and $y$ are orthogonal (or disjoint) as
elements of the {\it l}-group $G$, i.e. $|x| \wedge |y| = 1$,
where $|x| = x \vee x^{-1} = x_{+}x_{-}$. Consequently, $G$
is a $\perp$-group by \cite[Proposition 3.1.3.]{BKW}  
Notice that the $\perp$-group $G$ above is simplicial iff
it is Abelian, freely generated by the minimal positive
elements.

Moreover, in a {\it l}-group $G$, the following are satisfied:
$x \subset y \Lra x^{-1} \subset y^{-1}$, $x \cap y = 1 
\Lra xz \cap yz \subset z$, and $x \cup x^{-1} \neq \infty \Lra
x = 1$.

Inspired by the properties above satisfied by {\it l}-groups,
we introduce the following classes of median groups :

\begin{de} A median group $G$ is said to be an $A_i$-{\em group}, 
$i = 1, 2, 3, 4$, if $G$ satisfies the corresponding condition

$(A_1)$ $x \cup y \neq \infty$ and $x^{-1} \subset y^{-1} \Lra
x \subset y$

$(A_2)$ $x \cap y = x^{-1} \cap z = y^{-1} \cap z = 1 \Lra 
xz \cap yz \subset z$

$(A_3)$ $x \cup x^{-1} \neq \infty \Lra x^2 = 1$

$(A_4)$ $x \cup x^{-1} \neq \infty \Lra x = 1$
\smallskip

The median group $G$ is said to be an $A$-{\em group} if
$G$ is a $\perp$-group and also an $A_i$-group for $i = 1, 2, 4$.
\end{de}

Notice that the class of $A_i$-groups, $i = 1, 2, 3, 4$, as well as
the class of $A$-groups, is closed under arbitrary products.

\begin{rems} \em
(1) The $l$-groups are $A$-groups.

(2) Obviously, the locally linear median groups are $\perp$-groups.
They are also $A_2$-groups. Indeed, assume that $G$ is a locally 
linear median group, and $x, y, z \in G$ satisfy the identities
$x \cap y = x^{-1} \cap z = y^{- 1} \cap z = 1$. Setting 
$u := xz \cap yz$, we have by assumption $x, u \in [1, xz]$ and
$y, u \in [1, yz]$. As $G$ is locally linear, we distinguish
the following four cases.

$(i)\,$ : $u \subset x, u \subset y$. Then $u \subset x \cap y = 1$,
and hence $u = 1 \subset z$.

$(ii)\,$ : $x \subset u, y \subset u$. Then either $x \subset y$
or $y \subset x$, therefore either $x = 1$ or $y = 1$ since 
$x \cap y = 1$ by assumption. Consequently, $u \subset z$.

$(iii)\,$ : $x \subset u, u \subset y$. Thus $x \subset y$,
therefore $x = x \cap x \subset x \cap y = 1$, so $x = 1$, and
hence $u \subset z$.

$(iv)\,$ : $u \subset x, y \subset u$. As in $(iii)$, we get
$u \subset z$ as desired.

On the other hand, the locally linear $A_1$-groups are obviously
$A_3$-groups, but they are not necessarily $A_4$-groups;
for instance, the cyclic group of order $2$ satisfies 
$(A_i),\,i = 1, 2, 3$, while $(A_4)$ is not satisfied.

(3) The locally linear median groups, and hence the $\perp$-groups too,
are not necessarily $A_1$-groups. To provide an example of a 
locally linear median group which is not an $A_1$-group, we define
a semilattice operation $\cap$ on the set $\Z$ of integers as
follows :
\medskip

$n \cap m =
\left\{  
\begin{array}{lc}
\mbox{min} (n, m) & \mbox{if} \;\; n, m \;\; \mbox{are even} \geq 0\\
\mbox{max} (n, m) & \mbox{if} \;\; n, m \;\; \mbox{are either even} \leq 0
\;\; \mbox{or odd}\\
n & \mbox{if} \;\; n \geq 0 \;\; \mbox{is even, and} \;\; m \;\; 
\mbox{is odd}\\
m & \mbox{if} \;\; m \geq 0 \;\; \mbox{is even, and} \;\; n \;\; 
\mbox{is odd}\\
0 & \mbox{otherwise}
\end{array} \right. $
\medskip

It follows that $n \subset m$ iff one of the following four
conditions is satisfied :

(i) $0 \leq n \leq m$ and $n, m$ are even

(ii) $m \leq n \leq 0$ and $n, m$ are even

(iii) $m \leq n$ and $n, m$ are odd  

(iv) $n \geq 0$ is even and $m$ is odd.

One checks that $\Z$ with the usual addition and the operation
$\cap$ as defined above becomes a locally linear median group $G$ which is not
an $A_1$-group since $m \cup n = m \neq \infty, -m \subset -n$, and
$m \not \subset n$ whenever $m$ and $n$ are odd integers such that
$m < n$. Notice that the cell $[0, m] = \{n \in \Z\,|\,n \subset m\}$
is finite for $m$ even : $[0, m] = \{n \in 2\Z\,|\,0 \leq n \leq m\}$
for $m \geq 0$, resp. $[0, m] = \{n \in 2\Z\,|\,m \leq n \leq 0\}$
for $m \leq 0$, while $[0, m] = \{n \in 2\Z\,|\,n \geq 0\} \cup
\{n \in 2\Z + 1\,|\,n \geq m\}$ is infinite for $m$ odd. Since
$2 \cap - 2 = 0$ and $2n \subset 1$ for all $n \geq 0$, it follows
that $G$ is not Archimedean (cf. Definition 3.3.)
By contrast, the Archimedean $\perp$-groups are $A_1$-groups (see 
Proposition 3.4.)

(4) For a locally linear median group $G$, the following assertions
are equivalent.

(i) $G$ is an $A$-group.

(ii) $G$ is an $A_i$-group, $i = 1, 4$.

(iii) $\forall x \in G \sm \{1\}, x \not \subset x^{- 1}$, and
$\forall x \in G \sm \{1\}, y \in G, x \subset xy \Lra xy \not \subset y$.

(iv) $\forall x \in G \sm \{1\}, y, z \in G, [xy, xz] \not 
\sse [y, z]$.

(i) $\Llra$ (ii) follows by (2), (ii) $\Lra$ (iii) holds in all median
groups, while (iii) $\Lra$ (ii) holds in locally linear median groups.
On the other hand, (iii) $\Llra$ (iv) in locally linear median groups
by \cite[Lemma 3.2.]{Hyp}

(5) $(A_1)$ and $(A_2)$ do not imply $(\perp)$, resp. $(A_3)$.
For instance, let $G = \Z/4\Z$ be the cyclic group of order 4. The 
canonical order $\subset$ on $(G, S = \{1 {\mbox{mod}} 4\})$ makes $G$ 
a simplicial locally boolean $A_i$-group for $i = 1, 2$. However 
$G$ is not a $\perp$-group since $1 \perp -1$, while 
$1 \cup -1 = 2 \neq 0 = 1 + (-1)$. Moreover $G$ is not an 
$A_3$-group since $1 \cup -1 = 2 \neq \infty$ but $1 + 1 = 2 \neq 0$.
By contrast, the conditions $(A_i), i = 1, 2, 3$, are obviously
satisfied by any locally boolean $\perp$-group $G$ since, according
to \cite[Corollary 2.1.]{Rav}, $G$ is isomorphic to a subdirect product
of a power set $(\Z/2 \Z)^I$ with the canonical group and median 
operations.

(6) There exist $A_4$-groups (and hence $A_3$-groups too) which are not 
$\perp$-groups and $A_i$-groups for $i = 1, 2$. Indeed, given an 
ordered field $(K, \leq)$, let $G = K_{>0}$ denote the multiplicative 
group of positive elements of $K$, with the action  
$G \times K \longrightarrow K, (x, a) \mapsto xa$ of $G$ on 
the additive group $K$.

The total order on $K$ makes $G$ and $K$ locally linear $A_i$-groups 
for $i =1, 4$, and $G$ acts as a group of automorphisms of the 
median group $K$. The semidirect product $H := K \mathop {\rtimes} 
G$, with $(a, x) (b, y) := (a + xb, xy), (a, x) \cap(b, y) :=
(a \cap b, x \cap y)$, for $a, b \in K, x, y \in G$,  
is an $A_4$-group but it is not a $\perp$--group and an 
$A_i$-group for $i = 1, 2$. Indeed, assuming $(a, x) \in H$ such that 
$(a, x) \cup (a, x)^{-1} = (a, x) \cup (-x^{-1}a, x^{-1}) \neq 
\infty$ it follows that $x \cup x^{-1} \neq \infty$, therefore $x = 1$, and 
$a \cup -a \neq \infty$, and hence $a = 0$. Thus $(a, x) = (0, 1)$ is the 
neutral element of $H$, therefore $H$ is an $A_4$-group. To check that 
$H$ is not a $\perp$-group, let $a \in K, x \in G$ be such that 
$a \neq 0, x \neq 1$. Obviously, $(a, 1) \cap (0, x) = (0, 1)$ 
and $(a, 1) \cup (0, x) = (a, x) = (a, 1)(0, x) \neq (0, x)(a, 1) = 
(xa, x)$, so $H$ is not a $\perp$-group. To verify that $H$ is 
not an $A_1$--group, let $a \in K$ be such that $0 < a < 1$, 
and set $x := (1, 1), y := (a, a)$. We get $x \cup y = (1, a) \neq y$, 
though $x^{-1} = (-1, 1) \subset (-1, a^{-1}) = y^{-1}$. 
Finally, to check that $H$ is not an $A_2$-group, set 
$x := (2, 2^{-1}), y := (0,2), z := (1,1)$. We obtain 
$x \cap y = x^{-1} \cap z = y^{-1} \cap z = (0, 1)$ but 
$xz \cap yz = (2 + 2^{-1}, 2^{-1}) \cap (2, 2) = (2, 1) \not\subset z$.
\end{rems}

As we have seen in Remarks 2.3. $(5)$, $(A_2)$ does not imply 
$(\perp)$, however the converse is still open :
\smallskip

\noindent {\bf Question}. {\it Does the condition $(\perp)$ implies} $(A_2)$ ? 

Partial answers to the question above are provided by Remarks 2.3. 
$(1), (2), (5)$, and Corollary 3.5.

\begin{de} A $\perp$-group $G$ is called {\em Archimedean}
if for every $x \in G$ satisfying $\; x \cap x^{-1} = 1$,
i.e. $x \subset x^2$, and for every $y \in G$, 
there exists $n \geq 0$ such that $x^n \cap y = x^m \cap y$ 
for all $m \geq n$.
\end{de}

The Archimedean totally ordered groups, identified by H\" older's
theorem with subgroups of the additive ordered group $(\R, +)$ of 
reals, and the simplicial $\perp$-groups, i.e. the partially 
commutative Artin-Coxeter groups, are natural examples of 
Archimedean $\perp$-groups.

\begin{pr} Any Archimedean $\perp$-group is
an $A_i$-group for $i = 1, 3$.
\end{pr}

\bp Given an archimedean $\perp$--group $G$, let $x, y \in G$ be 
such that $x \cup y \neq \infty$ and $x^{-1} \subset y^{-1}$. To conclude 
that $G$ is an $A_1$-group we have to show that $x \subset y$. Setting 
$z := x \cap y, u := x^{-1}z, v := y^{-1}z$, it follows that 
$u^{-1} \perp v^{-1}$, therefore $u \perp v$ by \cite[Lemma 2.2.4.]{AC}

On the other hand, $u \subset x^{-1} \subset y^{-1}$ and 
$v \subset y^{-1}$ imply 
$u \cup v = u \bullet v = v \bullet u \subset y^{-1} = v \bullet z^{-1}$, 
therefore $u \subset z^{-1}$ and $u \cap u^{-1} \subset 
z^{-1} \cap u^{-1} = 1$, in particular $u^n \subset u^{n + 1}$ for 
all $n \geq 0$ by \cite[Lemma 2.2.3.]{AC} It remains to show 
by induction that $u^n \subset z^{-1}$ for all $n \geq 0$ 
to conclude thanks to the archimedeanity of $G$ that 
$u = 1$, i.e. $x = z \subset y$ as desired. Assuming 
$u^n \subset z^{-1} = u^n \bullet z'$ for some $n \geq 0$, 
we get $u^n \bullet u \bullet z' = 
u \bullet z^{-1} = x^{-1} \subset y^{-1} = v \bullet z^{-1} 
= v \bullet \mathop{\underbrace{u \bullet \ldots \bullet 
u}}\limits_{n \; \; {\mbox{factors}}} \bullet z' = 
u^n \bullet v \bullet z'$, 
therefore $u \subset z'$ and hence $u^{n+1} \subset z^{-1}$ as required.

To check that $G$ is an $A_3$-group, let $x \in G$ be such that 
$x \cup x^{-1} \neq \infty$, and let $y = x \cap x^{-1}$. Setting 
$u := x^{-1}y, v := xy$, we obtain $u^{-1} \perp v^{-1}$, therefore 
$u \perp v$ and $u \cup v = u \bullet v = v \bullet u$ by 
\cite[Lemma 2.2.4.]{AC} Thanks to the archimedeanity of $G$ it suffices 
to show that $(uv)^n \subset (uv)^{n+1} \subset y$ for 
all $n \geq 0$ to conclude that $v \perp u = v^{- 1}$, so $u = v = 1$,
and hence $x^2 = 1$ as desired. Since 
$x = y \bullet u^{-1} = v \bullet y^{-1}$ it follows by \cite[Lemma 2.1.]{Rav} 
and \cite[Lemma 2.2.4.]{AC} that $u \subset y$ and $v \subset y$, and hence 
$u \cup v = u \bullet v = v \bullet u \subset y$. 
Setting $y' := (uv)^{-1}y$, we get further $u \bullet y' = 
y'^{-1} \bullet v^{-1}$, therefore, again by \cite[Lemma 2.1.]{Rav}
and \cite[Lemma 2.2.4.]{AC}, 
$u \bullet v \subset y'^{-1}$. Setting $y'' := y'uv$, we obtain 
$y'' \bullet u^{-1} = v \bullet y''^{-1}$ and hence as above 
$u \bullet v \subset y''$, therefore 
$(uv)^2 = u \bullet v \bullet u \bullet v \subset y$. Thus by 
repeatedly applying the procedure above we obtain 
$(uv)^n \subset (uv)^{n+1} \subset y$ for all 
$n \geq 0$ as required. \ep

\begin{co} Any simplicial $\perp$-group is an 
$A_i$-group for $i = 1, 2, 3$.
\end{co}

\bp The cases $i = 1, 3$ are immediate by Proposition 3.4. 
since the simplicial $\perp$-groups are Archimedean. To prove 
the case $i = 2$, assume that $G$ is a simplicial 
$\perp$-group, and let $x, y, z \in G$ be such that 
$x \cap y = x^{-1} \cap z = y^{-1} \cap z = 1$. To show that 
$u:= (x \bullet z) \cap (y \bullet z) \subset z$, we argue by induction on 
the length $d := l(u)$ of $u$ over the generating set 
$\widetilde S = \{s \in G \setminus \{1\} \mid [1,s] = \{1,s\}\}$ of $G$. 
The case $d = 0$ is trivial, so let us assume $d \geq 1$, say 
$u = s \bullet v$ with $s \in \widetilde S$. We 
distinguish the following three possibilities :

$(i)\,$ : $s \subset x$, say $x = s \bullet x'$. As 
$s \subset y \bullet z$ 
and $s \cap y \subset x \cap y = 1$, we obtain 
$s \perp y$, and hence $s \subset z$, say $z = s \bullet z'$, 
and $y \bullet s = s \bullet y$. Simplifying with $s$, it follows 
that $v = x' \bullet s \bullet z' \cap y \bullet z'$. 
Notice that $x' \bullet s \cap y = 1$. Indeed, assuming the contrary, 
there exists $t \in \widetilde S$ such that 
$t \subset x' \bullet s \cap y$, and hence $t \subset x' \cap y$ and 
$s \perp t$ since $s \perp y$. Consequently, 
$t \subset x = s \bullet x'$, therefore $t \subset x \cap y = 1$, 
a contradiction. Since $l(v) = d - 1$ it follows by the induction 
hypothesis that $v \subset z'$, and hence 
$u = s \bullet v \subset s \bullet z' = z$.

$(ii)\,$ : $s \subset y$. We proceed as in the case $(i)$.

$(iii)\,$ : $s \cap x = s \cap y = 1$. It follows that $s \perp x, s \perp y$ 
and $s \subset z$, say $z = s \bullet z'$, therefore 
$v = x \bullet z' \cap y \bullet z'$. As $x \cap y = 1$ and 
$l(v) = d - 1$, the induction hypothesis implies $v \subset z'$ 
and hence $u \subset z$. \ep

\begin{co} The necessary and sufficient condition for 
a simplicial $\perp$-group to be an $A$-group is that it
is an $A_4$-group.
\end{co}

As an immediate consequence of \cite[Theorem 2.4.1.]{AC}, we obtain
the following characterization of right-angled Artin groups.

\begin{co} Given a group $G$ with a set 
$S \sse G$ of generators, let $\subset$ denote the partial order
on $G$ induced by the canonical length function on $(G, S)$. 
Then, the following assertions are equivalent.

$(1)$ $1 \not\in S, S \cap S^{-1} = \{s \in S \mid s^2 = 1 \}$, 
and the order $\subset$ makes $G$ an $A$-group.

$(2)$ $S \cap S^{-1} = \es$, and the order $\subset$ makes 
$G$ a $\perp$-group.

$(3)$ $(G, S)$ is a right-angled Artin group, i.e. $G$ 
admits the presentation 
$$G \cong <S; [s, t] = 1\,{\rm for}\, s, t \in S, s \neq t,
\, st = ts\,{\rm holds\, in}\, G>$$
\end{co}

Thus the right-angled Artin groups are identified 
with the simplicial $A$-groups.

We end the present section with four useful lemmas. 

\begin{lem} Let $G$ be an $A_1$-group, and let $x,y \in G$ be 
such that $xy = yx$. Then, $x^{-1} \cap y = 1 \Llra y^{-1} \cap x = 1$.
\end{lem}

\bp Let $x,y \in G$ be such that $xy = yx$ and $x^{-1} \cap y = 
1$, i.e. $xy = x \bullet y$. Setting $z = y^{-1} \cap x$, we get $yz = y 
\cap yx \subset yx = xy$, and $xyz = Y(x,xy,xyx) \subset xy$ since $x 
\cap xyx \subset x \subset xy$, therefore $yz \cup xyz \neq \infty$. As 
$x^{-1} \cap yz \subset x^{-1} \cap y = 1$, we get $(yz)^{-1} \subset 
(yz)^{-1}x^{-1} = (xyz)^{-1}$, and hence $yz \subset xyz \subset xy = yx$ by 
$(A_1)$. Consequently, $z \subset x^{-1}z \subset  x^{-1}$. By symmetry, 
it follows that $z \subset y$ too, therefore $z \subset x^{-1} \cap y = 
1$, i.e. $y^{-1} \cap x = 1$, as desired. \ep

\begin{lem} Assume that $G$ is a $\perp$-group satisfying $(A_1)$ 
and $(A_2)$, and let $x,y,z \in G$ be such that $xyz = x \bullet y \bullet z = z 
\bullet y \bullet x$ and $x \cap z = 1$. Then, $xy = yx, xz = zx$, and $yz = zy$.  
\end{lem}

\bp First note that $x \perp z$, and hence $xz = zx$ since $G$ is a $\perp$-group
by assumption. Consequently, $u:= x^{-1}yx = z^{-1}yz$, so we have to show that 
$u = y$. As $G$ is a $\perp$-group, $x \perp z$ implies $y \bullet x = x \bullet u$ and 
$y \bullet z =  z \bullet u$, and hence $y = y \bullet x \cap y \bullet z = 
x \bullet u \cap z \bullet u \subset u$ by $(A_2)$. On the other hand, $x 
\perp z$ implies $x^{-1} \perp z^{-1}$ by \cite[Lemma 2.2.4.]{AC}, therefore $u^{-1} = 
u^{-1} \bullet x^{-1} \cap u^{-1} \bullet z^{-1} = x^{-1} \bullet y^{-1} 
\cap z^{-1} \bullet y^{-1} \subset y^{-1}$ by $(A_2)$ again. As 
$y \subset u$ and $u^{-1} \subset y^{-1}$, it follows by $(A_1)$ that 
$u = y$ as required. \ep

\begin{lem} Assume that $G$ is a $\perp$-group satisfying 
$(A_1)$ and $(A_2)$. Then, for all $x \in G, [1,x] \cap Z_G(x) = \{y \in G 
\mid y \subset x, xy = yx\}$ is a sublattice of the cell $[1,x]$.
\end{lem}

\bp Let $y,z \in [1,x] \cap Z_G(x), u := y \cap z, y':= u^{-1}y$ 
and $z' := u^{-1}z$. As $y' \subset u^{-1}x, z' \subset u^{-1}x$ and $y' 
\cap z' = 1$, we get $y' \perp z'$, therefore $y \cup z = y 
\mathop{\vee}\limits_{x} z = u \bullet y' \bullet z' = u \bullet z' 
\bullet y'$ since $G$ is a $\perp$-group. Set $v := (y \cup z)^{-1}x$. 
Since $y$ and $z$ belong to $Z_G(x)$, and $G$ satisfies $(A_1)$, 
it follows by Lemma 3.8. that 
$z' \bullet v \bullet u \bullet y' = u \bullet y' \bullet z' \bullet 
v = x = u \bullet z' \bullet y' \bullet v = y' \bullet v \bullet u 
\bullet z'$. According to Lemma 3.9. we get $z'vu = vuz'$ and $ 
y^{\prime}vu=vuy^{\prime}$, therefore $x=u y^{\prime} z^{\prime} v= 
z^{\prime}vuy^{\prime} = z^{\prime} y^{\prime} vu = y^{\prime} z^{\prime} 
vu= y^{\prime} vuz^{\prime} = vuy^{\prime} z^{\prime} $, and hence $y 
\cap z$ and $y \cup z$ belong to $Z_G(x)$ as desired. \ep

\begin{lem} Let $G$ be an Abelian median group. Then, the following
 assertions are equivalent.

$(1)\,$ $G$ is an $A$-group.

$(2)\,$ $G$ is an $A_1$-group, and $x \cap x^{- 1} = 1$, i.e.
$x \subset x^2$, for all $x \in G$.
\end{lem}

\bp $(1) \Lra (2)$ : We have only to show that $x \cap x^{- 1} = 1$ for
all $x \in G$. Let $x \in G$, and set $u := x \cap x^{- 1}$. As $G$
is Abelian and satisfies $(A_1)$, it follows by Lemma 3.8. that
$u, u^{- 1} \subset x$, therefore $u \cup u^{- 1} \neq \infty$, and
hence $u = 1$ by $(A_4)$.

$(2) \Lra (1)$ : We have to show that $G$ is a $\perp$-group satisfying
$(A_2)$ and $(A_4)$. 

Assuming that $x \perp y$, let us show that 
$u := x \cup y = xy$. As $x, y \subset u$ and $G$ is an Abelian
$A_1$-group, it follows by Lemma 3.8. that 
$x^{- 1}, y^{- 1} \subset u^{- 1}$, therefore $x, y \subset xy = yx$
since $x^{- 1} \cap y, y^{- 1} \cap x \subset u \cap u^{- 1}$ and
$u \cap u^{- 1} = 1$ by assumption. Consequently, $u = Y(x, y, xy) =
xy (x^{- 1} \cap y^{- 1}) = xy$ since $x \perp y \Lra x \cap y = 1 \Lra
x^{- 1} \cap y^{- 1} = 1$ again by Lemma 3.8. Thus $G$ is a $\perp$-group.

To show that $G$ satisfies $(A_2)$, let $x, y, z \in G$ be such that
$x \cap y = x^{- 1} \cap z = y^{- 1} \cap z = 1$.  As $G$ is an Abelian
$A_1$-group, we deduce by Lemma 3.8. that $z \subset xz = zx$ and
$z \subset yz = zy$, and hence $xz \cap yz = z(x \cap y) = z \subset z$
as desired.

Finally, to show that $G$ satisfies $(A_4)$, let $x \in G$ be such
that $u := x \cup x^{- 1} \neq \infty$. As $G$ is a
$\perp$-group, we get $u = x (x \cap x^{- 1})^{- 1} x^{- 1}$,
and hence $u = 1$ as required, since $G$ is Abelian and 
$x \cap x^{- 1} = 1$ by assumption.
\ep

\begin{co} Let $G$ be an Abelian locally linear median group.
Then, the following assertions are equivalent.

$(1)\,$ $G$ is an $A$-group.

$(2)\,$ $x \cap x^{- 1} = 1$ for all $x \in G$.

$(3)\,$ There exist only two opposite total orders on $G$
making $G$ a totally ordered Abelian group whose associated
median group is the given locally linear median group $G$.
\end{co}

\bp
$(1) \Llra (2)$ : By Lemma 3.11., we have to show that 
$(2) \Lra (A_1)$ in Abelian locally linear median groups.
Let $x, y \in G$ be such that $x \cup y \neq \infty$ and
$x^{- 1} \subset y^{- 1}$. By locally linearity, $x \cup y \neq
\infty \Lra$ either $x \subset y$ or $y \subset x$. In the
former case we are done, while in the latter case, 
$y \subset x \Llra x^{- 1} y \subset x^{- 1} \Lra
x^{- 1} y \subset y^{- 1} \Lra x^{- 1} y \cap y \subset
y^{- 1} \cap y = 1 \Llra  y^{- 1} \subset y^{- 1} x^{- 1} y =
x^{- 1}$, therefore $x = y \subset y$ as desired.

For $(2) \Llra (3)$ see \cite[Corollary 3.4., Remark]{Hyp}
\ep 

\bigskip


\section{Cyclically reduced elements in $A$-groups}

$\quad$ In the rest of this Section, as well as in Sections 5 and 6, 
$G$ will denote an arbitrary $A$-group.

The basic notion of a cyclically reduced word in a free group
extends naturally to $A$-groups as follows.

\begin{de} An element $w \in G$ is said to be 
{\em cyclically reduced} if $w \cap w^{-1} = 1$.
\end{de}

Among the $A$-groups for which all elements are cyclically reduced,
we mention the $l$-groups and the Abelian $A$-groups (by Lemma 3.11.)

\begin{lem} Given $w \in G$, let $u := w \cap w^{-1}$ 
and $v := u^{-1}wu$. Then, $u$ is the unique element of 
$G$ for which $v$ is  cyclically reduced and 
$w = u \bullet v \bullet u^{-1}$.
\end{lem}

\bp By definition of $u$ we obtain 
$w = u \bullet (u^{-1}w) = (wu) \bullet u^{-1}$, in particular 
$wu \subset w$. To show that $u \subset wu$ set 
$u^{\prime} := u \cap wu$. It follows that 
$(u^{\prime -1}u) \bullet (u^{-1}w) = 
(u^{\prime -1}wu) \bullet u^{-1}$, therefore 
$u^{\prime-1}u \subset u^{-1}$ since $G$ is a $\perp$-group 
and $u^{\prime -1}u \cap u^{\prime -1}wu = 1$. 
As, on the other hand, 
$(u^{\prime -1}u)^{-1} = u^{-1}u^{\prime} \subset u^{-1}$ 
since $u^{\prime} \subset u$, we get 
$(u^{\prime -1}u) \cup (u^{\prime -1}u)^{-1} \neq \infty$, 
and hence $u = u^{\prime} \subset wu$ by $(A_4)$. To show 
that $v$ is cyclically reduced, set 
$v^{\prime} := v \cap v^{-1}$ and 
$w^{\prime} := v^{\prime -1}vv^{\prime}$. With the argument 
above we obtain $w = u \bullet v \bullet u^{-1} = 
u \bullet v^{\prime} \bullet w^{\prime} \bullet v^{\prime -1} 
\bullet u^{-1}$, therefore 
$u \subset u \bullet v^{\prime} \subset w \cap w^{-1} = u$, 
i.e. $v^{\prime} = 1$ as required.

To prove the uniqueness part of the statement, let $s \in G$ 
be such that $s \subset ws \subset w$ and $t: = s^{-1}ws$ 
is cyclically reduced. As $s \subset w \cap w^{-1}$, 
it remains to check that $ws \subset w^2$ to conclude that 
$s = w \cap w^{-1}$. As $s^{-1}w  \cap 
s^{-1}w^{-1} = t \bullet s^{-1} \cap t^{-1} \bullet s^{-1}$ 
and $t \cap t^{-1} = 1$, it follows by $(A_2)$ that 
$s^{-1}w \cap s^{-1}w^{-1} \subset s^{-1} \cap s^{-1}w = 1$, 
and hence $ws \subseteq w^2$ as desired. \ep


\begin{lem} For all $w \in G$ and for all natural numbers 
$n,m \geq 1, w^n \cap w^{-m} = w \cap w^{-1}$.
\end{lem}

\bp Let $u := w \cap w^{-1}$ and $v := u^{-1}wu$. 
By Lemma 4.2., $v \cap v^{-1} = 1$, therefore
$v^n \cap v^{-m} = 1$ for $n,m \in \N$
according to \cite[Lemma 2.2.3.]{AC}, in particular, $v^n$ is 
cyclically reduced for all $n \in \Z$. To conclude 
that $w^n \cap w^{-m} = u$ for all $n \geq 1, m \geq 1$, 
it suffices to show that $w^n = u \bullet v^n \bullet u^{-1}$ 
for all $n \geq 1$. Indeed the last condition implies the
identity $w^n \cap w^{-m} = u \bullet 
(v^n \bullet u^{-1} \cap v^{-m} \bullet u^{-1}) = u$ since 
$t:= v^n \bullet u^{-1} \cap v^{-m} \bullet u^{-1} \subset u^{-1}$ 
by $(A_2)$, whence 
$t = t \cap u^{-1} \subset v^n \bullet u^{-1} \cap u^{-1} = 1$. 
To check that $w^n = u \bullet v^n \bullet u^{-1}$ for 
$n \geq 1$, we argue by induction on $n$. The case $n = 1$ 
is assured by Lemma 4.2., so assuming 
$w^n = u \bullet v^n \bullet u^{-1}$ for some $n \geq 1$, we have to 
show that $w^{n+1} = u \bullet v^{n+1} \bullet u^{-1}$. As 
$u^{-1} \cap v^n = 1$ by the induction hypothesis and $G$ 
is a $\perp$-group, it follows that $u^{-1} \cap v^{n+1} = 
u^{-1} \cap v^n \bullet v \subset u^{-1} \cap v = 1$,
i.e. $u \subset uv^{n+1}$. Thus it remains to show that 
$s:= v^{-n-1} \bullet u^{-1} \cap u^{-1} = 1$. As 
$s \subset v^{-n-1} \bullet u^{-1}, 
v^{-1} \subset v^{-n-1} \bullet u^{-1}$ and 
$s \cap v^{-1} \subset u^{-1} \cap v^{-1} = 1$, it follows 
that $s \perp v^{-1}$, therefore $s \perp v$ by 
\cite[Lemma 2.2.4.]{AC} Consequently, $v \bullet u^{-1} = v \bullet s 
\bullet (s^{-1} u^{-1}) = s \bullet v \bullet (s^{-1} u^{-1})$, 
and hence $s \subset u^{-1} \cap v \bullet u^{-1} = 1$, 
i.e. $s = 1$ as required. \ep

\begin{co} For $w \in G$ and $0 \neq n \in \Z, w^n$ is cyclically 
reduced if and only if $w$ is cyclically reduced.
\end{co}

\begin{co} The $A$-groups are torsion-free.
\end{co}

\bp Let $w \in G$ and $n \geq 1$  be such that 
$w^n = 1$. By Lemma 4.3. we get 
$w \cap w^{-1} = w^n \cap w^{-1} = 1$, and hence 
$w \subset w^n = 1$ by \cite[Lemma 2.2.3.]{AC},
so $w = 1$. \ep

\begin{rem} \em
Let $(G,S)$ be a right-angled Artin group. By Corrolaries 3.7.
and 4.5., $G$ is torsion-free. For $w \in G, n \geq 1$,
it follows by Lemmas 4.2. and 4.3. that 
$l(w^n) = 2l(w \cap w^{- 1}) + nl((w \cap w^{- 1})^{- 1}w(w \cap w^{- 1}))$.
Given $w \in G$ and $s \in \widetilde S = S \cup S^{-1}$, set 
$u := w \cap w^{-1}, v := u^{-1}wu, w^{\prime} := sws^{-1}, 
u^{\prime} := w^{\prime} \cap 
w^{\prime -1}, v^{\prime} := u^{\prime -1}w^{\prime}u^{\prime}$. 
We distinguish the following four cases:

Case $(1)$ $sw = s \bullet w$ and 
$ws^{-1} = w \bullet s^{-1}$ : Then, by the condition 
$(F^{\prime})$ \cite[1.5.]{AC}, either 
$w^{\prime} = s \bullet w \bullet s^{-1}$ in which case 
$u^{\prime} = s \bullet u, v^{\prime} = v$, and 
$l(w^{\prime}) = l(w) + 2$ or $sw = ws$, 
i.e. $w^{\prime} = w$.

Case $(2)$ $s^{-1} \subset u = w \cap w^{-1}$ : Then, 
$u^{\prime} = su, v^{\prime} = v$, and 
$l(w^{\prime}) = l(w) - 2$.

Case $(3)$ $s^{-1} \subset w$ and $ws^{-1} = w \bullet s^{-1}$ : 
As $s^{-1} \subset w = u \bullet v \bullet u^{-1}$ 
and $s^{-1} \cap u \subset s^{-1} \cap w^{-1} = 1$, 
it follows that $s^{-1} \perp u$, therefore 
$s^{-1} \perp u^{-1}$ by \cite[Lemma 2.2.4.]{AC}, and hence 
$v = s^{-1} \bullet (sv)$ by $(\perp)$. Consequently 
$u^{\prime} = sus^{-1} = u, v^{\prime} = (sv) \bullet s^{-1}$, 
and $l(w^{\prime}) = l(w)$.

Case $(4)$ $ s^{-1} \subset w^{-1}$ and $sw = s \bullet w$ : 
Applying Case $(3)$ to $w^{-1}$ we obtain
$v = (vs^{-1}) \bullet s, u^{\prime} = sus^{-1} 
= u, v^{\prime} = s \bullet (vs^{-1})$, and $l(w^{\prime}) = l(w)$.

The discussion above implies that for given cyclically reduced 
elements $w, w^{\prime} \neq 1$ in $(G, S)$, the 
necessary and sufficient condition for $w$ and $w^{\prime}$ to be 
conjugate is that $l(w) = l(w^{\prime})$ and there exist a 
sequence $w_1 = w, w_2, \ldots , w_n = w^{\prime}$ of length 
$n \leq l(w)!$ and $s_i \in \widetilde S$ such that 
$s_i \subset w_i$ and $w_{i+1} = (s_i^{-1}w_i) \bullet s_i$ 
for $i<n$. Consequently, as the word problem on $(G, S)$ 
is solvable \cite{Tits}, \cite[1.5.]{AC}, it follows that 
the conjugacy problem on $(G, S)$ is solvable too.
\end{rem}

We end this section with four technical lemmas concerning
$A$-groups which will be useful later.

\begin{lem} Given the elements $x, w$ of the 
A-group $G$, set
$u := w \cap w^{-1}, v := u^{-1}wu, w^{\prime} := xwx^{-1}, u^{\prime} := 
w^{\prime} \cap w^{\prime -1}, v^{\prime} :=
u^{\prime -1}w^{\prime}u^{\prime}$, and 
$z := x^{-1} \cap w^{-1}x^{-1}$. If $xw = x \bullet w$ 
then $xw^2 = x \bullet w^2, z \perp u, u^{\prime} = (xz) \bullet u$, 
and $v^{\prime} = z^{-1}vz$.

In particular, if in addition $w$ and $w'$ are both 
cyclically reduced, then $z = x^{- 1}$, i.e. $x w = w' \bullet x$.
\end{lem}

\bp We proceed step by step as follows:

1) $xw^2 = x \bullet w^2$ : By Lemmas 4.2. and 4.3., 
$w = u \bullet v \bullet u^{-1}$ and 
$w^2 = u \bullet v \bullet v \bullet u^{-1}$. Setting 
$s:= x^{-1} \cap w^2$, we get 
$s \cap u \bullet v \subset x^{-1} \cap w = 1$, whence 
$s \perp u \bullet v$ and $s \subset v \bullet u^{-1}$ by 
$(\perp)$. On the other hand, $s \perp u \bullet v$ implies 
$s \perp v$ and $s \perp u^{-1}$ by \cite[Lemma 2.2.4.]{AC}, 
therefore $s = s \cap v \bullet u^{-1} = 1$, as desired.

2) $z \perp u$, in particular, $z \bullet u = u \bullet z$ : 
As $z \subset w^{-1} x^{-1}$ and $u \subset w^{- 1} 
\subset w^{- 1}x^{- 1}$, it follows that $z \cup u \neq \infty$.
Consequently, $z \perp u$ since $z \cap u \subset x^{-1} \cap w = 1$.

3) $xzu = (xz)\bullet u$ : As $x^{-1} \cap u \subset x^{-1} \cap w = 1$, 
we get $z \subset x^{-1} \cap z \bullet u = x^{-1} \cap u \bullet z 
\subset z$ by 2) and $(\perp)$, therefore $z^{-1}x^{-1} \cap u = 
z^{-1}Y(x^{-1},z,z \bullet u) = 1$, as required.

4) Set $a:= x^{-1} \cap w^{-1} = x^{-1} \cap u \bullet v^{-1} \bullet 
u^{-1}$. As $a \subset z$ and $z \perp u$ by 2), it follows that 
$a \perp u$, and hence $a \subset v^{-1}$ by \cite[Lemma 2.2.5.]{AC} Setting 
$x = y \bullet a^{-1}$ and $v = c \bullet a^{-1}$, we obtain

(4.1.) $y^{-1} \cap u \bullet c^{-1} \bullet u^{-1} = 1$.

Thus $z = a \bullet b$, where $b: = y^{-1} \cap u \bullet c^{- 1} \bullet 
u^{-1} \bullet x^{-1} \subset x^{-1}$ and 
$b \perp u \bullet c^{- 1} \bullet u^{-1}$ by (4.1.). As 
$z = a \bullet b \subset x^{-1}$ and $b \subset x^{-1}$, 
we get $b \subset a \bullet b$ by $(A_1)$, i.e. 
$z = a \bullet b = b \bullet a^{\prime}$, with 
$a^{\prime} = b^{-1}ab$. Setting $y = y^{\prime} \bullet b^{-1}$, 
we obtain $x = y^{\prime} \bullet b^{-1} \bullet a^{-1} = 
y^{\prime} \bullet a^{\prime -1} \bullet b^{-1}$ and

(4.2.) $y^{\prime -1} \cap u \bullet c^{-1} \bullet u^{-1} \bullet 
a^{\prime} \bullet y^{\prime -1} = 1$.

5) It remains to show that $u^{\prime} = xzu = xuz$, since 
then we get 
$v^{\prime} = u^{\prime -1}w^{\prime}u^{\prime} = z^{-1}vz$, 
as required. As $u^{\prime} = w^{\prime} \cap w^{\prime -1} = 
xY(wx^{-1}, x^{-1}, w^{-1}x^{-1})$, we have to show that 
$Y(wx^{-1}, x^{-1}, w^{-1}x^{-1}) = z \bullet u = 
u \bullet z = z \cup u$. Since $x^{-1} \cap w^{-1}x^{-1} = z$ 
by definition, it remains to check that 
$wx^{-1} \cap x^{-1} = b \subset z$ and 
$wx^{-1} \cap w^{-1}x^{-1} = u \bullet b \subset u \bullet z$.

5.1.) We deduce from (4.1.) and (4.2.) that 
$$wx^{-1} \cap x^{-1} = u \bullet c \bullet u^{-1} \bullet b 
\bullet y^{\prime -1} \cap b \bullet a^{\prime} \bullet y^{\prime -1} 
= b \bullet (u \bullet 
c \bullet u^{-1} \bullet y^{\prime -1} \cap a^{\prime} \bullet 
y^{\prime -1}) = b.$$ 

5.2.) It follows that $wx^{-1} \cap w^{-1}x^{-1} = u \bullet c \bullet 
u^{-1} \bullet b \bullet y^{\prime -1} \cap  u \bullet a \bullet c^{-1} 
\bullet u^{-1} \bullet b \bullet a^{\prime} \bullet y^{\prime -1} = u 
\bullet b \bullet t = b \bullet u \bullet t$, where 
$t = c \bullet u^{-1} \bullet y^{\prime -1} 
\cap a^{\prime} \bullet c^{-1} \bullet u^{-1} \bullet a^{\prime} \bullet 
y^{\prime -1}$, since $b \perp u \bullet 
c^{-1} \bullet u^{-1}$. By 1) we get $w^{-2}x^{-1} = w^{-2} \bullet 
x^{-1} = u \bullet v^{-2} \bullet u^{-1} \bullet x^{-1} = u \bullet a 
\bullet c^{-1} \bullet a \bullet c^{-1} \bullet u^{-1} \bullet b \bullet 
a^{\prime} \bullet y^{\prime -1} = u \bullet a \bullet b \bullet c^{-1} 
\bullet a^{\prime} 
\bullet u^{-1} \bullet a^{\prime} \bullet y^{\prime -1} = z \bullet u 
\bullet c^{-1} 
\bullet a^{\prime} \bullet c^{-1} \bullet u^{-1} \bullet a^{\prime} \bullet 
y^{\prime -1}$, therefore 
$c \bullet u^{-1} \cap a^{\prime} \bullet c^{-1} \bullet u^{-1} 
\bullet a^{\prime} \bullet y^{\prime -1} = 1$, and hence 
$t \subset y^{\prime -1} \cap u^{- 1} \bullet y'^{- 1}$. As $u^{- 1} \cap t = 1$,
it follows that $t \perp u^{\pm 1}$, in particular, $t \bullet u = u \bullet t$. 
On the other hand, since $z = b \bullet a' \perp u$, it follows by 
\cite[Lemma 2.2.5.]{AC} that $a' \perp u^{\pm 1}$, and hence we deduce
from (4.2) that $c \bullet u^{-1} \bullet y^{\prime -1} \cap a^{\prime} = 1$, 
whence $t \subset y^{\prime -1} \cap c^{-1} \bullet u^{-1} \bullet a^{\prime} 
\bullet y^{\prime -1}$. As $t \bullet u = u \bullet t$, we get 
$t \subset y^{\prime -1} \cap u \bullet 
c^{-1} \bullet u^{-1} \bullet a^{\prime} \bullet y^{\prime -1} = 1$ (by 
(4.2)), therefore $t = 1$ as desired. \ep

\begin{rem} \em With the notation above, we have $z^{-1} \bullet v = 
v^{\prime} \bullet z^{-1}, v = c \bullet a^{-1}$ and 
$v^{\prime} = a^{\prime -1} \bullet c$, where 
$a = x^{-1} \cap w^{-1}, a^{\prime} = z^{-1}az = b^{-1}ab, b = 
wx^{-1} \cap x^{-1}$. 
\end{rem}

\begin{lem} For all $x,y,w \in G, x \subset 
wx, y \subset wx$ and $y \cap w = 1 \Lra y \subset x$.
\end{lem}

\bp By Lemma 4.2., $x^{-1}wx = u \bullet v \bullet u^{-1}$, 
where $u = x^{-1}wx \cap x^{-1}w^{-1}x$ and $v = u^{-1} x^{-1}wxu$. 
Applying Lemma 4.7. to the reduced pair $(x,x^{-1}wx)$, we get 
$w = u^{\prime} \bullet v^{\prime} \bullet u^{\prime -1}$, 
where $u^{\prime} = w \cap w^{-1} = (xz) \bullet u, 
z = x^{-1} \cap x^{-1} w^{-1} = x^{-1}(x \cap w^{-1})$, so 
$xz = x \cap w^{-1}, v^{\prime} = u^{\prime -1} w u^{-1}$. 
Consequently, $xz \subset wxz \subset w$. 
As $y \subset wx = x \bullet (x^{-1}wx) = (wxz) \bullet z^{-1}$ 
and $y \cap wxz \subset y \cap w = 1$, it follows that 
$y \perp wxz$, in particular, $y \perp xz$ and $y \subset z^{-1}$. 
Since $x = (xz) \bullet z^{-1} = (xz) \bullet y \bullet (y^{-1}z^{-1}) 
= y \bullet (xz) \bullet (y^{-1} z^{-1})$, we get 
$y \subset x$ as desired.
\ep

\begin{lem} For $x,y,w \in G, x \subset w$ and 
$x \perp y \Lra x \subset ywy^{- 1}$.
\end{lem}

\bp Setting $ z = y^{-1} \cap w, z \subset y^{-1}$ and 
$x \perp y$ imply by \cite[Lemma 2.2.4.]{AC} $x \perp z$ and $x \perp yz$. 
As $x \subset w$ and $z \subset w$, we get 
$w = x \bullet z \bullet w^{\prime} = z \bullet x \bullet w^{\prime}$, 
and $yw = (yz) \bullet (x \bullet w^{\prime}) = x \bullet (yz) 
\bullet w^{\prime}$. It follows that $(yw)^{-1} \cap y^{-1} = 
w^{\prime -1} \bullet (yz)^{-1} \bullet x^{-1} \cap y^{-1} = 
w^{\prime -1} \bullet (yz)^{-1} \cap y^{-1}$ since 
$x^{-1} \perp y^{-1}$ by \cite[Lemma 2.2.4.]{AC} We conclude that 
$x \subset ywy^{-1}$ as required. \ep

\begin{lem} Given $x \in G$ and a cyclically 
reduced element $w \in G$, set 
$a := x \cap w, b := x \cap w^{-1}, y := a^{-1}b^{-1}x$, and 
$u := a^{-1}wb$. The necessary and sufficient condition for 
the conjugate $x^{-1}wx$ of $w$ to be cyclically reduced 
is that $y = ay \cap by$ and $y \perp u$.
\end{lem}

\bp As $w \cap w^{-1} = 1$, we obtain 
$a \perp b, x = a \bullet b \bullet y = b \bullet a \bullet y$, 
and $z:= a \bullet y \cap b \bullet y \subset y$ by $(A_2)$. 
Setting $ w = a \bullet w^{\prime}, b \subset w^{-1} = w^{\prime -1} 
\bullet a^{-1}$ and $b \perp a^{-1}$ (by \cite[Lemma 2.2.4.]{AC}) 
imply $b \subset w^{\prime -1}$, i.e. $w = a \bullet u \bullet b^{-1}$, 
and $a \bullet y \cap u^{-1} \bullet a^{-1} = 
b \bullet y \cap u \bullet b^{-1} = 1$, in particular, $z \cap u = 1$. 
We get $x \cap wx = a \bullet b \bullet y 
\cap a \bullet u \bullet a \bullet y = a \bullet s$, where 
$s:= b \bullet y \cap u \bullet a \bullet y \subset y$ by 
$(A_2)$ since $b \cap u \bullet a = 1$. By symmetry, we obtain 
$x \cap w^{-1}x = b \bullet t$, where 
$t:=a \bullet y \cap u^{-1} \bullet b \bullet y \subset y$. 
On the other hand, $wx \cap w^{-1}x = 
a \bullet u \bullet a \bullet y \cap b 
\bullet u^{-1} \bullet b \bullet y \subset y$ by $(A_2)$ 
again since 
$a \bullet u \bullet a \cap b \bullet u^{-1} \bullet b 
\subset w^2 \cap w^{-2} = 1$.

Assuming that $x^{-1}wx$ is cyclically reduced, it follows that 
$y^{-1}s, y^{-1}t \subset x^{-1}wx \cap x^{-1}w^{-1}x = 1$, 
therefore $y = s = t = z$, and $y \cap u = z \cap u = 1$. As 
$y = s \subset u \bullet a \bullet y$, we also get $y \perp u$.

Conversely, assuming $y = z$ and $y \bot u$, we get 
$wx \cap w^{-1}x = y \subset x, x \cap wx=a \bullet y$, 
and $x \cap w^{-1}x=b \bullet y$, and hence 
$Y(wx, x, w^{-1}x) = x \cap Y(wx, x, w^{-1}x)=(x \cap wx) \cup (x \cap 
w^{-1}x)=(ay)y^{-1}(by)=aby=x$, so $x \in [wx, w^{-1}x]$, i.e. 
$x^{-1}wx$ is cyclically reduced. 
\ep

\bigskip


\section{Preorders induced by elements of $A$-groups}

$\quad$ Given an element $w$ of the $A$-group $G$, let 
$\mathop{\preceq}\limits_w$ denote the binary relation 
defined by $x \mathop{\preceq}\limits_w y$ iff $y \in [x,wy]$, 
i.e. $x^{-1} y \subset x^{-1}wy$. Notice that 
$zx \mathop{\preceq}\limits_{zwz^{-1}} zy \Llra 
x \mathop{\preceq}\limits_w y$ for all $x, y, z \in G$.

\begin{lem} The relation $\mathop{\preceq}\limits_w $ 
is a preorder.
\end{lem}

\bp As the reflexivity is trivial, it remains to 
check the transitivity of $\mathop{\preceq}\limits_w$. 
Using a convenient conjugation, it suffices to show that 
$x \mathop{\preceq}\limits_w y$ whenever 
$x \mathop{\preceq}\limits_w 1$ (i.e. $x \cap w = 1)$ and 
$1 \mathop{\preceq}\limits_w y$ (i.e. $y \subset wy)$. 
By Lemma 4.9. we get $x \cap wy \subset y$, therefore 
$Y(x, y, wy) = y$, i.e. $x \mathop{\preceq}\limits_w y$
as required. \ep

Let $\mathop{\sim}\limits_w$ denote the equivalence 
relation induced by the preorder $\mathop{\preceq}\limits_w$.

\begin{lem} The necessary and sufficient 
condition for $x \mathop{\sim}\limits_w y$ is that 
$[x, wy]=[y, wx]$, i.e. $y^{-1} x \bot y^{-1} wy$.
\end{lem} 

\bp The non-trivial implication to prove is 
$x \mathop{\sim}\limits_w y \Lra
y^{-1}x \bot y^{-1}wy$. Without loss we may assume that 
$y=1$, so we have to show that $w \subset wx$ provided 
$x \subset wx$ and $x \cap w=1$. Applying Lemma 4.7. to 
the reduced pair $(x, x^{-1} wx)$, we get 
$w=(wxz) \bullet(xz)^{-1}$ and $wx = x \bullet 
(x^{-1} wx)=(wxz) \bullet z^{-1}$, where 
$z=x^{-1} \cap x^{-1}w^{-1}=x^{-1}(x \cap w^{-1})$. 
As $x \subset wx, wxz \subset wx$ and 
$x \cap wxz \subset x \cap w =1$, it follows by $(\perp)$ 
that $x \subset z^{-1}$, whence $x = z^{-1}$ by $(A_1)$ 
since $z \subset x^{-1}$. Consequently, $wx = w \bullet x$
as desired. \ep

\begin{lem} If $x \mathop{\preceq}\limits_{w} y$ 
and $z \in [x, y]$, then $x \mathop{\preceq}\limits_{w} z$ and 
$z \mathop{\preceq}\limits_{w} y$.
\end{lem}

\bp We may assume without loss that $x = 1$, 
so we have to show that $y \subset wy$ and $x \subset y$ 
imply $x \subset wz$ and $y \in [z, wy]$. As 
$z \subset y \subset wy$, we get obviously $y \in [z, wy]$, so 
it remains to prove that $z \subset wz$. In other words, 
we have to show that for $x, y, w \in G$, the pair 
$(x, ywy^{-1})$ is reduced whenever the triple $(x, y, w)$ 
is reduced, i.e. $x \subset xy \subset xyw$. Assuming 
that the triple $(x, y, w)$ is reduced, let us first 
apply Lemma 4.7. to the reduced pair $(y, w)$. Setting 
$u = w \cap w^{-1}, v = u^{-1}wu, z = y^{-1} \cap w^{-1} \bullet 
y^{-1}, u^{\prime} = w^{\prime} \cap w^{\prime -1}, v^{\prime} = 
u^{\prime -1}w^{\prime}u^{\prime}$, we obtain $z \perp u, u^{\prime} = 
(yz) \bullet u$ and $v^{\prime} \bullet z^{-1} = z^{-1} \bullet v$. 
Consequently, $1 = x^{-1} \cap y \bullet w = 
x^{-1} \cap (yz) \bullet z^{-1} \bullet u 
\bullet v \bullet u^{-1} = 
x^{-1} \cap u^{\prime} \bullet v^{\prime} \bullet 
u^{-1} \bullet z^{-1}$, therefore 
$s:= x^{-1} \cap ywy^{-1} = x^{-1} \cap u^{\prime} 
\bullet v^{\prime} \bullet u^{\prime -1} = x^{-1} \cap 
u^{\prime}
\bullet v^{\prime} \bullet u^{-1} \bullet (yz)^{-1} \subset (yz)^{-1}$ 
by $(\perp)$. As $s$ and $yz \subset u^{\prime}$ are bounded 
above by $u^{\prime} \bullet v^{\prime} \bullet u^{-1} \bullet (yz)^{-1}$, 
and $s \cap yx \subset x^{-1} \cap u^{\prime} = 1$, we get 
$s \perp yz$, and hence $s \perp (yz)^{-1}$ by \cite[Lemma 2.2.4.]{AC} 
Since, on the other hand, $s \subset (yz)^{-1}$, it follows 
that $s = 1$, i.e. the pair $(x, ywy^{-1})$ is 
reduced as desired. \ep

The preorder $\mathop{\preceq}\limits_{w}$ is 
compatible with the arboreal structure of $G$, as follows :

\begin{pr} Given $x, y, a, b \in G, Y(a, b, y) 
\mathop{\preceq}\limits_{w} Y(a, b, x)$ whenever 
$y \mathop{\preceq}\limits_{w} x$. In particular, the 
equivalence relation $\mathop{\sim}\limits_{w}$ is a 
congruence on the underlying median set of $G$.
\end{pr}

\bp We may assume that $y = 1$, so we have to show that 
$x \subset wx \Lra Y(a, b, x) \in [a \cap b, Y(wa, wb, wx)]$. 
Setting $c:= Y(a, b, x) = (a \cap b) \cup (c \cap x)$, and 
$d := (a \cap b) \cap (c \cap x) = a \cap b \cap x$, we 
get $c = (a \cap b)d^{-1} (c \cap x) = (c \cap x)d^{-1}(a \cap b), 
d^{-1}x \subset d^{-1}wx$, and 
$d^{-1} (a \cap b) \perp d^{-1}(c \cap x)$, whence 
$(a \cap b)^{-1}d \perp d^{-1}(c \cap x)$ by \cite[Lemma 2.2.4.]{AC} 
As $d^{-1}(c \cap x) \subset d^{-1}x$, it follows by 
Lemma 5.3. that $d^{-1}(c \cap x) 
\subset d^{-1}w(c \cap x)$, therefore
$d^{-1}(c \cap x) \subset 
((a \cap b)^{-1}d) (d^{-1}w(c \cap x)) 
(d^{-1}(a \cap b)) 
= (a \cap b)^{-1}wc$, according to Lemma 4.10. 
Multiplying with $a \cap b$, we obtain $c \in [a \cap b, wc]$
as required. \ep

For any  $x \in G$, let $\widetilde x^w$ denote the 
$\mathop{\sim}\limits_{w}$ - class of $x$. Note that 
$\widetilde x^w$ is a convex subset of $G$. For 
$x, y \in G$ such that $x \mathop{\preceq}\limits_{w} y$, 
set $<x, y>_{w}:= \{z \in G \mid x \mathop{\preceq}\limits_{w} z\; 
{\mbox{and}}\; z \mathop{\preceq}\limits_{w} y\}$.

\begin{co} Let $x, y \in G$ be such that 
$x \mathop{\preceq} \limits_{w} y$. The convex subset 
$<x, y>_w$ is the disjoint union 
$\mathop{\bigsqcup}\limits_{z \in [x,y]} \widetilde z^w$.
\end{co}

\begin{co} $\widetilde 1^w = 
\{x \in G \mid x \perp w\}$ 
is a convex subgroup of $G$.
\end{co}

\bp The closure of $\widetilde 1^w$ under the 
operation $x \mapsto x^{-1}$ follows by \cite[Lemma 2.2.4.]{AC} 
To check the closure of $\widetilde 1^w$ under multiplication, 
let $x,y \in \widetilde 1^w$, and set $z := x^{-1} \cap y$. 
Thus $xy = (xz) \bullet (z^{-1}y)$, where $xz$ and $z^{-1}y$ 
belong to $\widetilde 1^w$ by \cite[Lemma 2.2.5.]{AC}, so we may 
assume from the beginning that $xy = x \bullet y$, with 
$x, y \in \widetilde 1^w$. As $x \perp w$ and $y \perp w$ 
we obtain $w \cap x \bullet y = w^{-1} \cap x \bullet y = 1$ 
by $(\perp)$, therefore $xy \perp w$ since 
$w \bullet x \bullet y = x \bullet w \bullet y = 
x \bullet y \bullet w$. \ep

\begin{lem}For all $n \geq 1$, the preorder 
$\mathop{\preceq}\limits_{w^n}$ induced by $w^n$ equals 
$\mathop{\preceq}\limits_{w}$.
\end{lem}

\bp It suffices to show that 
$x \mathop{\preceq}\limits_{w} 1 \Llra
x \mathop{\preceq}\limits_{w^n} 1$, i.e. $x \cap w = 1 \Llra
x \cap w^n = 1$. It remains to argue as in the step 1) 
of the proof of Lemma 4.7. \ep

\bigskip


\section{Foldings induced by elements of $A$-groups}

$\quad$ For any element $w$ of an $A$-group $G$, let us define the mapping 
$\varphi_w : G \lra G$ by $\varphi_w(x) = Y(wx, x, w^{-1}x)$. 
The main goal of this section is to show that $\varphi_w$ 
is a folding of the underlying median set of $G$. In the particular
case of locally linear $A$-groups, this fact is a consequence of
\cite[Proposition 2.7.]{Hyp}  

Notice that $\varphi_{xwx^{-1}}(y) = x \varphi_w (x^{-1}y)$ 
for all $x,y \in G$.

\begin{lem} $X_w := \{x \in G\,|\,\varphi_w(x) = x\} = 
\{x \in G\,|\,x \mathop{\preceq}\limits_{w} wx\}$ 
is a convex subset of $G$.
\end{lem}

\bp Notice that $X_w$ consists of those $x \in G$ 
for which the conjugate $x^{-1}wx$ of $w$ is cyclically reduced. 
Assuming without loss that $w$ is cyclically reduced, we have 
to show that $y^{-1}wy$ is cyclically reduced whenever 
$x^{-1}wx$ is cyclically reduced and $y \subset x$. Setting 
$a = x \cap w, b = x \cap w^{-1}, x^{\prime} = a^{-1}b^{-1}x$, 
and $u = a^{-1}wb$, it follows by Lemma 4.11. that 
$a \bullet x^{\prime} \cap b \bullet x^{\prime} = x^{\prime}$ 
and $x^{\prime} \perp u$. Setting 
$c = y \cap w = y \cap a, d = y \cap w^{-1} = y \cap b, y^{\prime} 
= c^{-1} d^{-1}y, \alpha = c^{-1}a, \beta = d^{-1}b$, 
we get $c \perp d, \alpha \perp d, \beta \perp c$ 
(by \cite[Lemma 2.2.5.]{AC} since $a \perp b$), and 
$x = a \bullet b \bullet x^{\prime} 
= c \bullet \alpha \bullet d \bullet \beta \bullet x^{\prime} 
= c \bullet d \bullet \alpha \bullet \beta \bullet x^{\prime}$, 
so $y^{\prime} \subset \alpha \bullet \beta \bullet x^{\prime}$. 
As $d \bullet y^{\prime} \cap \alpha = 1$ and  $\alpha \perp d$, 
it follows that $y^{\prime} \perp \alpha$, and by symmetry, 
$y^{\prime} \perp \beta$, so $y^{\prime} \subset x^{\prime}$ 
by $(\perp)$. Setting $z = y^{\prime -1}x^{\prime}$, we get 
$y^{\prime} \subset x^{\prime} \subset a \bullet x^{\prime} 
= c \bullet \alpha \bullet y^{\prime} \bullet z = 
c \bullet y^{\prime} \bullet \alpha \bullet z$, therefore 
$y^{\prime} \subset c \bullet y^{\prime}$ by $(A_1)$, 
and also, by symmetry, $y^{\prime} \subset d \bullet y^{\prime}$. 
As $c \perp d$ we obtain 
$c \bullet y^{\prime} \cap d \bullet y^{\prime} = y^{\prime}$ 
by $(A_2)$. On the other hand, $y^{\prime} \perp u$ since 
$x^{\prime} \perp u$ and $y^{\prime} \subset x^{\prime}$. 
As we already know that $y^{\prime} \perp \alpha$ and 
$y^{\prime} \perp \beta$ it follows by Corollary 5.6. that 
$y^{\prime} \perp \alpha \bullet u \bullet \beta^{-1} (=c^{-1}wd)$, 
therefore $y^{-1}wy$ is cyclically reduced according 
to Lemma 4.11. \ep

\begin{lem} $\varphi_w(G) = X_w$. 
\end{lem}

\bp The inclusion $X_w \subseteq \varphi_w(G)$ is trivial. 
To prove the opposite inclusion, let $x \in G$, and set 
$w^{\prime} := x^{-1}wx$. By Lemma 4.3., $w^{\prime 2} 
\cap w^{\prime -2} = w^{\prime} \cap w^{\prime -1}$, 
therefore $\varphi_w(\varphi_w(x)) =$ 
$$Y(w \varphi_w(x), \varphi_w(x), w^{-1} \varphi_w(x)) 
= x Y(w^{\prime 2} \cap w^{\prime}, w^{\prime} 
\cap w^{\prime -1}, w^{\prime -1} \cap w^{\prime -2}) 
= x(w^{\prime} \cap w^{\prime -1}) = \varphi_w(x),$$
i.e. $\varphi_w(x) \in X_w$ for all $x \in G$. \ep

\begin{lem} If $w$ is cyclically reduced, then 
$wx \cap w^{-1}x \subset x$ for all $x \in G$.
\end{lem}

\bp Setting $a = x \cap w, b = x \cap w^{-1}, 
y = a^{-1}b^{-1}x$ and $u = a^{-1} wb$, we get as in 
the initial part of the proof of Lemma 4.11. that 
$z:= wx \cap w^{-1}x = a \bullet u \bullet a \bullet y \cap b 
\bullet u^{-1} \bullet b \bullet y \subset y$, so 
$z = a \bullet u \bullet a \bullet z \cap b \bullet u^{-1} 
\bullet b \bullet z$ by $(A_2)$. As $ a \bullet z \cap b^{-1} = 1$ 
we get $z \cap w = a \bullet u 
\bullet a \bullet z \cap b \bullet u^{-1} \bullet b \bullet z \cap 
a \bullet u \bullet b^{-1} = a \bullet u \cap b \bullet 
u^{-1} \bullet b \bullet z$, and hence (according to 
\cite[Lemma 2.2.5.]{AC}) $z \cap w = z \cap a$, since 
$a \bullet u \cap b \bullet u^{-1} \bullet b \subset w \cap w^{-2} = 1$ 
and $z \cap a \bullet u \subset b \bullet u^{-1} \bullet b \bullet z 
\Lra (z \cap a \bullet u) \perp b \bullet u^{-1} \bullet b$.
By symmetry, it follows that $z \cap w^{-1} = z \cap b$. 
Setting $c = z \cap w, d = z \cap w^{-1}, \alpha = 
c^{-1}a, \beta = d^{-1}b, z^{\prime} = c^{-1}d^{-1}z$, we get 
$z = c \bullet d \bullet z^{\prime} = d \bullet c \bullet z^{\prime}$, 
and also $c \perp u$ and $d \perp u$ by \cite[Lemma 2.2.5.]{AC} 
As $z = wx \cap w^{-1}x \subset \varphi_w(x)$ it follows 
by Lemmas 6.1. and 6.2. that $z \in X_w$, i.e. $z^{-1}wz$ 
is cyclically reduced, and hence 
$c \bullet z^{\prime} \cap d \bullet z^{\prime} = z^{\prime}$ 
and $z^{\prime} \perp \alpha \bullet u \bullet \beta^{-1}$ 
according to Lemma 4.11. Setting $c^{\prime} = 
z^{\prime -1}cz^{\prime}, d^{\prime} = z^{\prime -1}dz^{\prime}$ 
and $y^{\prime} = z^{-1}y$, we get $x = a \bullet b 
\bullet y = a \bullet b \bullet z \bullet y^{\prime} 
= c \bullet \alpha \bullet d 
\bullet \beta \bullet c \bullet d \bullet z^{\prime} \bullet y^{\prime} 
= c \bullet d 
\bullet \alpha \bullet \beta \bullet z^{\prime} \bullet c^{\prime} 
\bullet d^{\prime} \bullet y^{\prime} 
= c \bullet d \bullet z^{\prime} \bullet \alpha \bullet \beta \bullet 
c^{\prime} \bullet d^{\prime} \bullet y^{\prime}$, therefore
$z = c \bullet d \bullet z^{\prime} \subset x$ 
as required. \ep

\begin{pr} $\varphi_w$ is a folding of $G$.
\end{pr}

\bp We may assume without loss that $w$ is 
cyclically reduced, i.e. $1 \in X_w$. By Lemmas 6.1. and 6.2., 
$X_w = \varphi_w(G)$ is a convex subset of $G$, so it remains to show 
that $[1,x] \cap X_w = [1, \varphi_w(x)]$ for all $x \in G$. 
As $\varphi_w(x) = (wx \cap w^{-1}x) \cup (x \cap \varphi_w(x))$, 
it follows by Lemma 6.3. that $\varphi_w(x) \subset x$, 
i.e. $\varphi_w(x) \in [1,x]$. Conversely, let 
$y \in [1,x] \cap X_w$. As $y^{-1}wy$ is cyclically reduced, 
it follows as above that 
$y^{-1} \varphi_w(x) = \varphi_{y^{-1}wy}(y^{-1}x) \subset y^{-1}x$, 
therefore $x^{-1} \varphi_w(x) \subset x^{-1}y \subset x^{-1}$, 
so $y \subset \varphi_w(x)$, i.e. 
$y \in [1, \varphi_w(x)]$ as desired. \ep

\begin{lem} $X_w$ is closed under the congruence 
$\mathop{\sim}\limits_{w}$. In particular, the folding 
$\varphi_w$ induces a folding of the quotient median set 
$G/\mathop{\sim}\limits_{w}$.
\end{lem}

\bp Let $x \in X_w$ and $y \in G$ be such that 
$x \mathop{\sim}\limits_{w}y$, i.e. $[x, wy] = [y, wx]$. 
Assuming $y \not\in X_w$, i.e. $y \not\in [wy, w^{-1}y]$, 
it follows by \cite[Corollary 5.2.2.]{Dual} that there is 
a prime convex subset $P$ of $G$ such that 
$wy \in P, w^{-1}y \in P$ and $y \not\in P$. As 
$[x, wy] = [y, wx]$ and $[x, w^{-1}y] = [y, w^{-1}x]$, 
it follows that  $x \not\in P, wx \in P$ and $w^{-1}x \in P$, 
therefore $x \not\in [wx, w^{-1}x]$, contrary to the assumption 
$x \in X_w$. \ep

Let $\mathop{\equiv}\limits_{w}$ denote the negation 
of the congruence $\mathop{\sim}\limits_{w}$ in the 
Heyting algebra $Cong\; (G)$, cf. 2.1., and let 
$\mathop{\ll}\limits_{w}$ denote the order on $G$ 
defined by $x \mathop{\ll}\limits_{w} y$ iff 
$x \mathop{\preceq}\limits_{w} y$ 
and $x \mathop{\equiv}\limits_{w} y$.

\begin{lem} On $X_w = X_{w^{-1}}$ 
the preorders $\mathop{\preceq}\limits_{w}$ and  
$\mathop{\preceq}\limits_{w^{-1}}$, as 
well as the orders $\mathop{\ll}\limits_{w}$ and 
$\mathop{\ll}\limits_{w^{-1}}$, are opposite.
\end{lem}

\bp We may assume without loss that $1 \in X_w$ 
and $x \in X_w$ such that $1 \mathop{\preceq}\limits_{w} x$, 
i.e. $x \subset wx$. Let $a = x \cap w, b = x \cap w^{-1}$ 
and $y = a^{-1}b^{-1}x$. Since $w$ and $x^{-1}wx$ are 
cyclically reduced, it follows as in the proof of Lemma 4.11. 
that $a \bullet b \bullet y =x=x \cap wx = a \bullet y$, 
therefore $b=x \cap w^{-1} = 1$, i.e. 
$x \mathop{\preceq}\limits_{w^{-1}} 1$ as required. \ep

\begin{lem} For all $x \in X_w, x \mathop{\ll}\limits_{w} wx$.
\end{lem}

\bp As $x \mathop{\preceq}\limits_w wx$ provided 
$x \in X_w$, it remains to show that $x \in X_w \Lra
x \mathop{\equiv}\limits_w wx$.
Let $y, z \in [x, wx]$ be such that 
$y \mathop{\sim}\limits_w z$, i.e. $[y, wz]=[z, wy]$. 
Assuming that $y \neq z$, it follows by \cite[Corollary 5.2.2.]{Dual}
that there is a prime convex subset $P$ of 
$G$ such that $y \in P$ and $z \not\in P$, therefore $wy \in P$ 
and $w z \not\in P$. On the other hand, as $y, z \in [x,wx]$, 
we distinguish the following two cases :

Case $(1)$ : $x \in P$ and $wx \not\in P$. As 
$wy \in [wx, w^2x] \cap P$, we get $w^2 x \in P$, therefore 
$wx \in [x, w^2 x] \subseteq P$, i.e. a contradiction.

Case $(2)$ : $x \not\in P$ and $w x \in P$.  As 
$wz \in [wx, w^2x] \backslash P$, we get $w^2 z \not\in P$, 
and hence $wx \in [x, w^2 x] \subseteq G \backslash P$, 
again a contradiction.

Consequently, $y = z$ as desired. \ep

\begin{lem} $\varphi_{w^n}= \varphi_w$ for all $n \neq 0$.
\end{lem}

\bp As $\varphi_w$ and $\varphi_{w^n}$ are foldings 
of $G$ it suffices to show that they have a common image, 
i.e. $X_w=X_{w^n}$. The equality above is now immediate by 
Corollary 4.4. \ep

\begin{lem} For all
$x \in G, x \mathop{\ll}\limits_{w} \varphi_w (x)$. 
In particular, the orders $\mathop{\ll}\limits_w$,  
$\mathop{\ll}\limits_{w^{-1}}$ and 
$\mathop{\leq}\limits_{\varphi_w (x)}$ coincide on the cell 
$[x, \varphi_w (x)]$.
\end{lem}

\bp First let us show that 
$x \mathop{\preceq}\limits_w \varphi_w (x)$, i.e. 
$\varphi_w(x) \in [x,w \varphi_w (x)] = [x, \varphi_w (wx)]$. 
Assuming the contrary, it follows by \cite[Corollary 5.2.2.]{Dual}
that there is a prime convex subset $P$ of 
$G$ such that $x \in P, \varphi_w (wx) \in P$, and 
$\varphi_w (x) \not\in P$. Consequently, $wx \not\in P$ and 
$w^2x \in P$, whence 
$\varphi_w (x) = \varphi_{w^2}(x)=Y(w^2x, x, w^{-2}x) \in P$, 
i.e. a contradiction. 

Next let us show that 
$x \mathop{\equiv}\limits_w \varphi_w (x)$. Assuming that there are 
$y, z \in [x, \varphi_w (x)]$ such that $y \neq z$ and y 
$\mathop{\sim}\limits_w z$, i.e. $[y, wz]=[z,wy]$, it 
follows by \cite[Corollary 5.2.2.]{Dual} again that there is 
a prime convex subset $P$ of $G$ such that 
$y \in P, wy \in P, z \not\in P$, and $wz \not\in P$, therefore 
$\varphi_w (y) \in P$ while $\varphi_w (z) \not\in P$, contrary 
to $\varphi_w(y) = \varphi_w (z) = \varphi_w(x)$ since 
$y, z \in [x, \varphi_w (x)]$ by assumption. \ep

\begin{rems}  \em
$(1)$ Being a folding, $\varphi_w$ induces according to 
2.2. a quasidirection 
$\mathop{\bullet}\limits_{\varphi_w}$ defined by 
$x \mathop{\bullet}\limits_{\varphi_w} y = 
Y(x, y, \varphi_w (x))$, whose associated preorder 
$\mathop{\preceq}\limits_{\varphi_w}$ 
is given by $x \mathop{\preceq}\limits_{\varphi_w} y$ 
iff $y \mathop{\bullet}\limits_{\varphi_w} x=y$ iff 
$y \in [x, \varphi_w(y)]$. Notice that 
$\mathop{\preceq}\limits_{\varphi_w}$ is finer than the 
preorders $\mathop{\preceq}\limits_w$ and 
$\mathop{\preceq}\limits_{w^{-1}}$.

$(2)$ Obviously, for all $w \in G$, the centralizer 
$Z_G(w)=\{x \in G | xw =wx\}$ of $w$ in $G$ is 
contained in the stabilizer $Stab \; (\varphi_w) = 
\{x \in G|\varphi_{xwx^{-1}}
= \varphi_w\} = \{x \in G|xX_w= X_w\}$ of 
$\varphi_w$ under the action from the left of $G$. 
However the converse is not necessarily true, as
for instance in non-commutative $l$-groups. Indeed,
if $G$ is an $l$-group, then $\varphi_w = \varphi_1 = 1_G$
for all $w \in G$.

$(3)$ Given a median set $X$, one assigns to any automorphism 
$s$ of $X$ the mapping $\varphi_s : X \lra X$, defined by
$\varphi_s(x) = Y(s(x), x, s^{- 1}(x))$. According to
Proposition 6.4., $\varphi_s$ is a folding whenever $X$ is
the underlying median set of an $A$-group $G$ and 
$s$ is the left translation $x \mapsto s(x) := wx$ by 
some element $w \in G$. An analogous situation is
provided by \cite[Theorem 6.6.]{AB} for a 
$\Lam$-{\em tree} (cf. \cite{MS}) $X$, where $\Lam$ is 
a totally ordered Abelian
group, and a {\em hyperbolic automorphism} $s$ of $X$.
In this case, the $s$-{\em axis} $X_s := \varphi_s(X)$ is 
identified with a convex subset of $\Lam$, and 
$s|_{X_s}$ is equivalent to a translation $x \mapsto
x + l(s)$, where $0 < l(s) = {\rm min}_{x \in X}\, d(x, s(x))$ 
( $= d(p, s(p))$ for some (for all) $p \in X_s$) is the 
{\em hyperbolic length} of $s$.

More generally, we can consider a {\em faithfully full
$\Lam$-metric median set} (cf. \cite[1.3.]{RS}), where 
$\Lam$ is an Abelian $l$-group, and an automorphism $s$ 
of $X$. Then $\varphi := \varphi_s$ is an endomorphism of the underlying
median set of $X$ satisfying the following equivalent conditions :

$(i)\, \varphi^3 = \varphi$ and $\varphi^2$ is a folding;

$(ii)\, \varphi^3 = \varphi$ and $\varphi(X)$ is a convex 
subset of $X$;

$(iii)\, \forall x, y, z \in X, \varphi(Y(x, y, \varphi(z))) = 
Y(\varphi(x), \varphi(y), z)$. 

The following assertions hold :

$(a)\, X_s := \varphi(X) = \varphi^2(X) =
\{x \in X\,|\,\varphi^2(x) = x\}$ is a retractible convex
subset of $X$;

$(b)\, \varphi|_{X_s}$ is an involutive automorphism of $X_s$,
and $\varphi|_{X_s} = s|_{X_s} \Llra {\rm Fix}(s^2) :=
\{x \in X\,|\,s^2(x) = x\} \neq \es$;

$(c)\, l(s) := (d(x, s^2(x)) - d(x, s(x))_+ \in \Lam_+$ does not depend
on the element $x \in X$, and $l(s) = 0 \Llra\,{\rm Fix}(s^2) \neq \es$;

$(d)\, [x, s(x)] = [\varphi(x), \varphi(s(x)) = s(\varphi(x))]$ 
for all $x \in X_s$;

$(e)\, d(x, \varphi(s(x))) = l(s)$ for all $x \in X_s$;

$(f)\,\forall x \in X, d(x, s(x)) = d(\varphi(x), \varphi(s(x))) +
2 d(x, \varphi^2(x)) = l(s) + d(\varphi(x), \varphi^2(x)) +
2 d(x, \varphi^2(x))$.

Details will be given in a forthcoming paper. Notice
also that the pair $(\varphi, \varphi)$, with $\varphi = \varphi_s$
as above, is a particular case of the so called {\em compatible pairs}
cf. \cite[Section 11]{UC2} which are a basic ingredient for the construction
of {\em universal coverings} relative to {\em median groupoids of median
sets} and {\em simplicial median groupoids of groups} 
\cite[Proposition 11.4., Theorem 14.1.]{UC2}, which extend
the {\em universal covering relative to a connected graph of groups}
\cite[Ch. I, Theorem 12]{Ser}
\end{rems}

\bigskip


\section{Quasidirections induced by elements of right-angled Artin groups}

$\quad$ In the rest of the paper we assume that $(G, S)$ is a right-angled 
Artin group. By Corollary 3.7., the partial order $\subset $  
induced by the canonical length function on $(G, S)$ makes $G$
a simplicial $A$-group, so we can apply 
to this special case the general theory developed in the 
previous sections.

\begin{lem} For all $w, x \in G$, there exists 
$y \in G$ such that $x \cap wy =y$.
\end{lem}

\bp We argue by induction on the length $d := l(x)$. 
The case $d = 0$ is trivial. Assuming that the equality 
$x \cap w y = y$ is satisfied for some $y \in G$, let 
$s \in \widetilde{S} = S \cup S^{-1}$ be such that 
$xs = x \bullet s$. If $x s  \cap wy = y$, then we have 
nothing to prove, so let us assume that 
$y \mathop{\subset}\limits_{\neq} xs \cap wy$, 
therefore, by $(\perp), s \bot y^{-1}x, x \bullet s
= y \bullet s \bullet (y^{-1}x), wy = y \bullet s 
\bullet (s^{-1}y^{-1} wy)$, and $xs \cap wy = ys$. We 
distinguish two cases :

Case $(1) : wys = (wy)\bullet s = 
y \bullet s \bullet (s^{-1}y^{-1} wy) 
\bullet s$. As $y^{-1}x \cap s^{-1}y^{-1}wy = 1$ 
and $y^{-1}x \bot s$, we get 
$y^{-1} x \cap (s^{-1} y^{-1} wy) \bullet s = 1$, 
and hence $xs \cap wys = ys$ as required.

Case $(2) : wys \subset wy$. Thus 
$s \subset  y^{-1} w^{-1} = 
(y^{-1}w^{-1}ys) \bullet s^{-1} \bullet y^{-1}$, therefore 
$s \subset y^{-1} w^{-1}ys$, since otherwise 
$s \subset s^{-1} \bullet y^{-1}$ by $(\perp)$, 
contrary to $(A_4)$. Consequently, 
$x s\cap wys = x s \cap wy = ys$ as desired. \ep

\begin{pr} For all $w \in G$, the preorder 
$\mathop{\preceq}\limits_{ w}$, defined by 
$x \mathop{\preceq}\limits_{w} y \Llra y \in [x, wy]$, 
determines a quasidirection $\mathop{\bullet}\limits_{w}$ 
on $G$.
\end{pr}

\bp By Proposition 5.4., the preorder
 $\mathop{\preceq}\limits_{w}$ is compatible with the 
arboreal structure on $G$. Moreover each pair $(x,y)$ 
of elements in $G$ is bounded above with respect to 
$\mathop{\preceq}\limits_{w}$. Indeed, by Lemma 7.1 
applied to the elements $x^{-1} wx$ and $x^{-1} y$, there is 
$z \in G$ such that $x^{-1}y \cap x^{-1}wxz = z$, 
i.e. $Y(x, y, wxz) = xz$, therefore 
$x \mathop{\preceq}\limits_{w} xz$ and 
$y \mathop{\preceq}\limits_{w} xz$.

To conclude that the binary operation 
$\mathop{\bullet}\limits_{w}$, defined by 
$a \mathop{\bullet}\limits_{w} b = 
\mathop{\vee}\limits_a U_{a,b}$ 
with $U_{a,b} = \{x \in [a,b]| 
a \mathop{\preceq}\limits_{w} x \; \mbox{and} \; b 
\mathop{\preceq}\limits_{w}x\}$, is a quasidirection on $G$, 
it suffices to show, according to Lemma 2.4., that 
$c \mathop{\preceq}\limits_{w} b \Lra 
a \mathop{\preceq}\limits_{w} b$ whenever the elements
$a, b, c \in G$ satisfy $[a, b] = \{a, c, b\}, 
c \not\in \{a,b\}$. Setting 
$s := c^{-1}a, t := c^{-1}b$, and $w^{\prime} := c^{-1}wc$, 
it follows by assumption that $s, t \in \widetilde{S}, s \cap t =1$, 
i.e. $s \neq t$, and moreover $s \cup t = \infty$, so either 
$t =s^{-1} $ or $st \neq ts$. Assuming that
$c \mathop{\preceq}\limits_{w} b$, i.e. $t \subset w^{\prime} t$, 
it follows that $s \cap w^{\prime} t=1$, as $s \cup t = \infty $, 
and hence $t \in [s, w^{\prime} t]$, i.e. 
$a \mathop{\preceq}\limits_{w}b$, as required. \ep

\begin{pr} Given $a \in X_w := \varphi_w(G)$, 
let $X_{w, a}$ denote the convex closure of the subset 
$\{w^n a |n \in \Z\}$, and let $\Psi_{w,a}$ denote the 
folding of $G$ associated to $X_{w, a}$. Then, the following
assertions hold.

$(1)\,X_{w, a} = \mathop{\bigcup}\limits_{n \geq 0} [ w^{-n} a, w^{n} a]$ 
is an unbounded distributive lattice with respect to 
the order $\mathop{\ll}\limits_w$, the join 
$x \mathop{\bullet}\limits_{w} y =
y \mathop{\bullet}\limits_{w} x$, and the 
meet $x \mathop{\bullet}\limits_{w^{-1}} y = y 
\mathop{\bullet}\limits_{w^{-1}}x$ for $x, y \in X_{w,a}$.

$(2)\,$ For all $x \in G, \Psi_{w,a} (x) = 
\displaystyle\lim_{n \rightarrow 
\infty } Y(w^{-n}a, x, w^n a)$, i.e. there exists 
$m \geq 0$ such that $\Psi_{w,a} (x) = Y(w^{-n} a, x, w^n a)$ 
for all $n \geq m$.

$(3)\, X_w$ is the closure of $X_{w,a}$ under the congruence 
$\mathop{\sim}\limits_{w}$.

$(4)\, X_{w,a}=X_w \cap \mathop{a}\limits^{\equiv_w }$.
\end{pr}

\bp $(1)\,$ Since $w^na \mathop{\ll}\limits_{w} w^m a$ for 
$n,m \in \Z, n \leq m$, by Lemma 6.7., and the orders 
$\mathop{\ll}\limits_{w}$ and $\mathop{\leq}\limits_{y}$  
coincide on $[x,y]=\{z \in G| x \mathop{\ll}\limits_{w} z 
\mathop{\ll}\limits_{w}y\}$ provided $x \mathop{\ll}\limits_{w} y$, 
it follows that $X_{w,a}$ is the union of the ascending chain of 
cells $[w^{-n}a, w^n a]$ for $n \geq 0$, which is directed by 
the order $\mathop{\ll}\limits_{w}$. As the orders 
$\mathop{\ll}\limits_{w}$ and $\mathop{\ll}\limits_{w^{-1}}$ 
are opposite on $X_w$ according to Lemma 6.6., we obtain the 
desired structure of distributive lattice on $X_{w, a}$.

$(2)\,$ Thanks to the definition of $\Psi_{w, a}$ and to $(1)$, for all 
$x \in G$, there is $m \geq 0$ such that 
$[a, \Psi_{w, a} (x)]=[a, x] \cap X_{w, a} 
=[a,x] \cap [w^{-n}a, w^n a]=[a, Y(w^{-n}a, x, w^n a)]$ 
for all $n \geq m$, therefore 
$\Psi_{w, a} (x)= \displaystyle\lim_{n \rightarrow \infty } 
Y(w^{-n}a, x, w^n a)$.

$(3)\,$ We have to show that for all $x \in X_w, x 
\mathop{\sim}\limits_{w} \Psi_{w, a} (x)$, i.e. 
$w^{-n}a \mathop{\preceq}\limits_{w} x 
\mathop{\preceq}\limits_{w} w^n a$ for large enough $n$. 
Let $m \geq 0$ be such that 
$\Psi_{w, a}(x) = Y(w^{-n}a, x, w^n a)$ for $n \geq m$. 
Assuming that $x \mathop{\not\preceq}\limits_{w} w^m a$, i.e. 
$w^m a \not\in [x, w^{m+1} a]$, it follows by \cite[Corollary 5.2.2.]{Dual}
that there is a prime convex subset 
$P$ of $G$ such that $[x,w^{m+1} a] \subseteq P$ while 
$w^m a \not\in P$, therefore $\Psi_{w,a} (x) = 
Y(w^{-m-1} a, x, w^{m+1} a) \in P$. As we also have 
$\Psi_{w, a} (x) = Y(w^{-m}a,x,w^ma) \in P$, and 
$w^m a \not \in P$, it follows that $w^{-m} a \in P$, 
and hence $w^m a \in [w^{-m}a, w^{m+1}a] \subseteq P$, 
i.e. a contradiction. On the other hand, as $x \in X_w$, 
we may interchange the roles of $a$ and $x$ to get some 
$k \geq 0$ subject to $a \mathop{\preceq}\limits_{w} w^kx$, i.e. 
$w^{-k} a \mathop{\preceq}\limits_{w}  x$. Taking $n=\max(m,k)$, 
it follows that $w^{-n} a \mathop{\preceq}\limits_{w} x 
\mathop{\preceq}\limits_{w} w^n a$ as required.

$(4)\,$ The inclusion $X_{w,a} \subseteq X_w \cap \stackrel{\equiv}{a^w}$ 
is immediate by Lemma 6.7. Conversely, assuming that 
$x \in X_w \cap \stackrel{\equiv}{a^w}$, it follows by $(3)$ 
that there is $y \in X_{w,a}$ such that 
$x \mathop{\sim}\limits_{w} y$, 
therefore $x = y \in X_{w, a}$ since we also have 
$x \mathop{\equiv}\limits_{w} a \mathop{\equiv}\limits_{w}y$. \ep

\begin{co} Let $a \in X_{w}$. Then, $X_w$ is the convex 
closure of the $Z_G(w)$-orbit of $a$.
\end{co}

\bp Obviously $u a \in X_w$ for all $u \in Z_G (w)$. 
If $x \in X_w$, then Proposition 7.3. provides a natural 
number $n$ and some $y \in [w^{-n}a, w^na]$ such that 
$x \mathop{\sim}\limits_{w} y$, i.e. $x^{-1} y \bot x^{-1}wx$ 
(by Lemma 5.2.), therefore $xy^{-1} \in Z_G(w)$. 
It follows that $x \in [ua, va]$ with $u = xy^{-1} w^{-n}, 
v = xy^{-1} w^n \in Z_G(w)$. \ep

\begin{co} Given $w \in G$ and $a \in X_w$, 
let $\varphi_{w, a}$, resp. $\Psi_{w, a}$, denote 
the folding of $G$ associated to the convex subset 
$\widetilde a^w$, resp. $X_{w, a}$. Then,
the median set morphism 
$X_w \longrightarrow \widetilde a^w \times X_{w, a}, x 
\mapsto (\varphi_{w, a}(x), \Psi_{w, a} (x))$ 
is an isomorphism, whose inverse sends a pair 
$(y, z) \in \widetilde a^w \times X_{w,a}$ to 
$y a^{-1} z= za^{-1}y$.
\end{co}

\bp By Proposition 7.3., $X_{w, x} = X_w \cap 
\mathop{x}\limits^{\equiv_w}$ for all $x \in X_w$, and  
$X_{w,x} \cap \widetilde y^w \neq \es$
for all $x, y \in X_w$. Thus we may apply Lemma 2.1. and 
Corollary 2.2. to conclude that the mapping above is an 
isomorphism of median sets, whose inverse sends a pair 
$(y, z) \in \widetilde a^w \times X_{w, a}$ to 
$\Psi_{w, y} (z) = \varphi_{w, z} (y)$. As 
$[y, z]= [\varphi_{w, y}(z), \Psi_{w, y} (z)] =
[a, \Psi_{w, y} (z)]$ by Lemma 2.1., we get 
$\Psi_{w, y}(z)=ya^{-1}z=z a^{-1} y$ as required. \ep

The next statement provides a description of the 
quasidirection $\mathop{\bullet}\limits_{w}$ by 
means of the folding $\varphi_w$.

\begin{co} Let $w \in G$. Then, $x \mathop{\bullet}\limits_{w} y = 
\displaystyle\lim_{n \rightarrow \infty } Y(x,y, w^n \varphi_w(x))$ 
for all $x, y \in G$.
\end{co}

\bp Let $x, y \in G$. As $x \mathop{\ll}\limits_{w} w^m 
\varphi_w (x)$ for all $m \geq 0$ by Lemmas 6.7 and 6.9., 
taking into account the definition of 
$\mathop{\bullet}\limits_{w}$, 
it suffices to show that 
$y \mathop{\preceq}\limits_{w} w^n \varphi_n (x)$ for large 
enough $n$. By Lemma 6.9., $y \mathop{\ll}\limits_{w} \varphi_w(x)$, 
so, in particular, $y \mathop{\preceq}\limits_{w} \varphi_w (y)$. 
On the other hand, according to Proposition 7.3. there exist 
$m \geq 0$ and $z \in [w^{-m} \varphi_w(x), w^m \varphi_w(x)]$ 
such that $\varphi_w(y) \mathop{\sim}\limits_{w} z$, whence 
$\varphi_w (y) \mathop{\preceq}\limits_{w} w^n \varphi_w(x)$ 
for all $n \geq m$, as required. \ep

\begin{co} For all $w \in G$, the folding 
$\varphi_w$, interpreted as a quasidirection through the 
convex embedding 
$Fold (G) \longrightarrow Dir(Fold(G)) \cong Qdir (G)$, 
is the join of the quasidirections $\mathop{\bullet}\limits_{w} $ 
and $\mathop{\bullet}\limits_{w^{-1}}$ in the directed 
median set $Qdir (G)$.
\end{co}

\bp By definition, the join of the 
quasidirections $\mathop{\bullet}\limits_{w} $ and
$\mathop{\bullet}\limits_{w^{-1}}$ is the quasidirection $\bullet 
:=Y(\mathop{\bullet}\limits_{w} , \mathop{\bullet}\limits_{w^{-1}}, 
\mathop{\bullet}\limits_{1})$ defined by $x \bullet y=Y(x 
\mathop{\bullet}\limits_{w} y, x \mathop{\bullet}\limits_{w^{-1}}y, x)$ 
for $x, y \in G$. Given $x, y \in G$, it follows by Corollary 7.6. 
that there is $n \geq 0$ such that $x \mathop{\bullet}\limits_{w} y
= Y(x, y, w^n\varphi_wx))$ and 
$x \mathop{\bullet}\limits_{w^{-1}} y = 
Y( x, y, w^{-n} \varphi_w(x))$. As 
$Y(w^n \varphi_w(x), w^{-n} \varphi_w(x), x) = 
\varphi_w(x)$ by Lemmas 6.7. and 6.9., we get $x \bullet y 
= Y(x, y, Y(w^n \varphi_w(x), w^{-n} \varphi_w(x), x))
= Y(x, y, \varphi_w(x))$, i.e. $\bullet $ is the quasidirection 
induced by the folding $\varphi_w$. \ep

\begin{co} Let $w \in G$. Then, the convex subset of 
$Fold \; (G)$ obtained by intersecting the cell 
$[\mathop{\bullet}\limits_{w}, \mathop{\bullet}\limits_{w^{-1}}]$ of 
$Qdir (G)$ with the convex subset 
$Fold (G)$ of $Qdir (G)$ consists of those foldings 
$\eta $ of $G$ for which $\eta(x) \in X_{w,\varphi_{w}(x)}$ 
for all $x \in G$.
\end{co}

\bp By definition, the intersection 
$[\mathop{\bullet}\limits_{w},
\mathop{\bullet}\limits_{w^{-1}}]\cap \; Fold \; (G)$ 
consists of the foldings $\eta $ of $G$ subject to 
$Y(x, y, \eta (x)) \in [x \mathop{\bullet}\limits_w y, x 
\mathop{\bullet}\limits_{w^{-1}} y]$ for all $x, y \in G$. 
Given such a folding $\eta$ and taking 
$x \in \eta (G)$ and $y = \varphi_w (x)$, we get 
$x = \varphi_w (x)$ since 
$x \mathop{\bullet}\limits_{w} \varphi_w(x) = x 
\mathop{\bullet}\limits_{w^{-1}}\varphi_w(x) = 
\varphi_w (x)$ by Lemma 6.9. Thus $\eta (G) \subseteq X_w$. 
On the other hand, taking $y = \eta(x)$, we obtain 
$\eta(x) \in [x \mathop{\bullet}\limits_{w} \eta (x), x
\mathop{\bullet}\limits_{w^{-1}} \eta (x)]$, therefore, by applying 
$\varphi_w$, we get $\eta (x)=\varphi_w(\eta(x)) \in [\varphi_w(x) 
\mathop{\bullet}\limits_{w} \eta(x), \varphi_w (x) 
\mathop{\bullet}\limits_{w^{-1}} \eta(x)] \subseteq X_w \cap 
\mathop{\varphi_w(x)}\limits^{\equiv_w}= X_{w,\varphi_w (x)}$ 
(by Proposition 7.3.$(4)$). Conversely, if the folding 
$\eta $ of $G$ satisfies the condition 
$\eta(x) \in X_{w, \varphi_w (x)}$ for all $x \in G$, 
then, thanks to Proposition 7.3.$(1)$, for each $x \in G$ 
there exists $m \geq 0$ such that 
$\eta (x) \in [w^{-n} \varphi_w(x), w^n \varphi_w(x)]$ for 
all $n \geq m$, therefore, by Corollary 7.6., $Y(x,y,\eta(x)) \in [x 
\mathop{\bullet}\limits_w y, x \mathop{\bullet}\limits_{w^{-1}}y]$ 
for all $y \in G$. \ep

\begin{pr} For $w, a \in G$, let 
$\mathop{\vee}\limits_{w; \; a}$ denote the direction on $G$ 
obtained by applying the folding of $Dir (G)$ induced by 
the quasidirection $\mathop{\bullet}\limits_{w}$ to the 
internal direction $\mathop{\vee}\limits_a$ on $G$ associated 
to $a$, and let $\mathop{\leq}\limits_{w; \; a}$ denote 
the associated order on $G$. Then, the following assertions hold.

$(1)$ For $x,y \in G, x \mathop{\vee}\limits_{w; a}y = 
\displaystyle\lim_{n \rightarrow \infty}
Y(x, y, Y(w^n \varphi_w(x), a, w^n \varphi_w(y)))$. In 
particular, the directions $\mathop{\vee}\limits_{w;a}$ and 
$\mathop{\vee}\limits_{w;\varphi_w (a)}$ coincide.

$(2)$ The ray from $a$ in the direction 
$\mathop{\vee}\limits_{w;a}$, namely 
$[a, \mathop{\vee}\limits_{w;a}): =
[a, \mathop{\vee}\limits_{w;a}] \cap 
G = \{x \in G\,|\,a \mathop{\leq}\limits_{w;a} x\}$, consists 
of those $x \in G$ for which $a \mathop{\ll}\limits_w x$;
the orders 
$\mathop{\leq}\limits_{w; a}, \mathop{\ll}\limits_{w; a},$ 
and the opposite of $\mathop{\leq}\limits_a$ coincide on 
$[a, \mathop{\vee}\limits_{w;a})$ 
making it a distributive lattice with the meet 
$\mathop{\vee}\limits_a$, the join $\mathop{\bullet}\limits_w$, 
and the least element $a$.

$(3)$ The mapping $G \rightarrow Dir (G), \; a \mapsto 
\mathop{\vee}\limits_{w: a}$ is a morphism of median sets 
inducing a convex embedding of $G/\equiv_w\,\, \cong\, X_w/\equiv_w $ 
into $Dir (G)$.

$(4)$ The quasidirection induced by the folding 
$\Psi_{w, \varphi_w (a)}$ of $G$ associated to 
the convex subset $X_{w, \varphi_w (a)}$ is the 
join in the directed median set $Qdir (G) $ 
of the directions $\mathop{\vee}\limits_{w; \; a}$ and 
$\mathop{\vee}\limits_{w^{-1}; \; a}$.

$(5)\,$ $(\mathop{\vee}\limits_{w; \; a}, 
\mathop{\vee}\limits_{w^{-1}; \; a}):= 
[\mathop{\vee}\limits_{w; \; a}, 
\mathop{\vee}\limits_{w^{-1}; \; a}] \cap G = 
X_{w, \varphi_w (a)}$.
\end{pr}

\bp $(1)$ By definition, $x \mathop{\vee}\limits_{w; a} y 
= (x \mathop{\bullet}\limits_w y) \mathop{\vee}\limits_a (y 
\mathop{\bullet}\limits_w x)$ for all $x, y \in G$, and 
hence the equality stated in $(1)$ is immediate by Corollary 7.6. 
As $\varphi_w $ is a folding, we get 
$\mathop{\vee}\limits_{w; a} = \mathop{\vee}\limits_{w; 
\varphi_w (a)}$.

$(2)$ As we also have $x \mathop{\vee}\limits_{w; a} y= (x 
\mathop{\vee}\limits_a y) \mathop{\bullet}\limits_w x 
\mathop{\bullet}\limits_w y$, it follows that 
$a \mathop{\vee}\limits_{w; a} x= a \mathop{\bullet}\limits_w x$, 
therefore $a \mathop{\leq}\limits_{w; a} x \Llra 
a \mathop{\ll}\limits_w x$, as desired. For 
$x,y \in [a, \mathop{\vee}\limits_{w;a})$ we get 
$(x \mathop{\vee}\limits_a y)  
\mathop{\bullet}\limits_w x = 
Y(x, y, a \mathop{\bullet}\limits_w x) = 
Y(x, y, x) =x$, whence $x \mathop{\vee}\limits_{w; a} y = x 
\mathop{\bullet}\limits_w y = y \mathop{\bullet}\limits_w x$. 
In particular, for $x, y \in [a, \mathop{\vee}\limits_{w;a}), x 
\mathop{\leq}\limits_{w,a} y \Llra x \mathop{\ll}\limits_w y
\Llra x \in [a,y] \Llra y \mathop{\leq}\limits_a x$.

$(3)$ The compatibility with the median operations on $G$ 
and $Dir (G)$ of the mapping 
$a \mapsto \mathop{\vee}\limits_{w;a}$ is obvious from the 
definition of $\mathop{\vee}\limits_{w;a}$. Moreover, for 
$a, b \in  G$ and $D \in Dir (G)$, denoting by 
$\vee $ the direction 
$Y(\mathop{\vee}\limits_{w; a}, \mathop{\vee}\limits_{w; b}, D)$, 
we get $x \vee y = (x \mathop{\vee}\limits_{w; a} y) 
\mathop{\vee}\limits_D (x \mathop{\vee}\limits_{w; b} y) = 
Y(x \mathop{\bullet}\limits_w y, y  
\mathop{\bullet}\limits_w x, a \mathop{\vee}\limits_D b)$ 
for all $x,y \in G$, therefore
$\vee = \mathop{\vee }\limits_{w; (a \mathop{\vee}\limits_D b)}$. 
Thus the image of the morphism $a \mapsto 
\mathop{\vee}\limits_{w; a}$ is a convex subset of $Dir (G)$. 
It remains to show that $\mathop{\vee}\limits_{w; a} = 
\mathop{\vee}\limits_{w; b} \Llra a \mathop{\equiv}\limits_w b$. 
Assuming that $\mathop{\vee}\limits_{w; a}= \mathop{\vee}\limits_{w; b}$, 
let $c := a \mathop{\vee}\limits_{w; a} b$. By $(2)$ we get 
$ a \mathop{\ll}\limits_w c$ and $ b \mathop{\ll}\limits_w c$, 
whence $a \mathop{\equiv}\limits_w c \mathop{\equiv}\limits_w b$. 
Conversely, as $a  \mathop{\equiv}\limits_w b \Llra a 
\mathop{\bullet}\limits_w b = b \mathop{\bullet}\limits_w a, a 
\mathop{\equiv}\limits_w b$ implies $a \mathop{\bullet}\limits_w  
b \in [a, \mathop{\vee}\limits_{w; a}) \cap [b, 
\mathop{\vee}\limits_{w; b})$ by $(2)$, therefore, by
$(2)$ again, $[a \mathop{\bullet}\limits_{w}  b, 
\mathop{\vee}\limits_{w; a}) = \{x \in G\,|\, a 
\mathop{\bullet}\limits_w  b \mathop{\ll}\limits_w x\} = 
[a \mathop{\bullet}\limits_w  b, \mathop{\vee}\limits_{w; b})$, 
i.e. $\mathop{\vee}\limits_{w; a} = \mathop{\vee}\limits_{w; 
b}$ as required.

$(5)$ By $(1)$ we may assume that $a \in X_w$. For any $x \in 
(\mathop{\vee}\limits_{w; a}, \mathop{\vee}\limits_{w^{-1}; a})$ 
we get $x = a \mathop{\vee}\limits_x x \in [a \mathop{\vee}
\limits_{w; a} x, a \mathop{\vee}\limits_{w^{-1}; a}x] 
\subseteq X_{w, a}$ since by $(1)$ there 
is $n \geq 0$ such that $a \mathop{\vee}\limits_{w; a} x = 
Y(a, x, Y(w^n a, a, w^n \varphi_w (x))) \in [a, w^n a]$, 
and similarly 
$a \mathop{\vee}\limits_{w^{-1}; a} x \in [a, w^{-n} a]$. 
Conversely, as $\mathop{\vee}\limits_{w; b} = 
\mathop{\vee}\limits_{w; a}$ for all 
$b \in X_{w, a} \subseteq \mathop{a}\limits^{\equiv_w}$ 
by $(3)$, it suffices to show that 
$a \in [\mathop{\vee}\limits_{w; a}, 
\mathop{\vee}\limits_{w^{-1}; a}]$, i.e. 
$x \mathop{\vee}\limits_a y \in [x \mathop{\vee}\limits_{w; a}y, 
x \mathop{\vee}\limits_{w^{-1}; ay}]$ 
for all $x, y \in G$. For $x,y \in G$, let 
$u := \mathop{\vee}\limits_{w; a} \{x, y,a\}$ and 
$v := \mathop{\vee}\limits_{w^{-1}; a} \{x, y, a\}$. It follows 
as above that $u$ and $v$ belong to $X_{w, a}$, and 
$a \in [v, u]$. Consequently, 
$x \mathop{\vee}\limits_a y = Y(x, y, u 
\mathop{\vee}\limits_a v) = (x \mathop{\vee}\limits_u y) 
\mathop{\vee}\limits_a(x \mathop{\vee}\limits_v y) = (x 
\mathop{\vee}\limits_{w; a} y) \mathop{\vee}\limits_a (x 
\mathop{\vee}\limits_{w^{-1}; a} y) \in [x \mathop{\vee}
\limits_{w; a} y, x \mathop{\vee}\limits_{w^{-1}; a} y]$
as desired.

$(4)$ We may assume that $a \in X_w$ by $(1)$. By definition we have 
to show that $Y(\mathop{\vee}\limits_{w; a}, x, 
\mathop{\vee}\limits_{w^{-1}; a}) = 
\Psi_{w, a}$ ($x$) for all $x \in G$, i.e. 
$\Psi_{w, a} (x) \in [ \mathop{\vee}\limits_{w; a}, 
\mathop{\vee}\limits_{w^{-1}; a}] \cap 
[x, \mathop{\vee}\limits_{w; a}] \cap [x, 
\mathop{\vee}\limits_{w^{-1}; a}]$. By $(5)$ we get 
$\Psi_{w, a} (x) \in [\mathop{\vee}\limits_{w; a}, 
\mathop{\vee}\limits_{w^{-1}; a}]$, while 
$\Psi_{w, a} (x) \in [a, x] \cap [a 
\mathop{\vee}\limits_{w; a} x, x]$ (since $a 
\mathop{\vee}\limits_{w; a} x \in X_{w, a})$ implies 
$Y(x, \mathop{\vee}\limits_{w; a}, \Psi_{w, a} (x)) = 
Y(x, \mathop{\vee}\limits_{w; a}, Y(a, x, \Psi_{w, a} (x))) = 
Y(a \mathop{\vee}\limits_{w; a} x, x, \Psi_{w, a} (x)) = 
\Psi_{w, a} (x)$, i.e. 
$\Psi_{w, a} (x) \in [x, \mathop{\vee}\limits_{w; a}]$. 
Similarly, we obtain $\Psi_{w, a} (x) \in [x, 
\mathop{\vee}\limits_{w^{-1}; a}]$. \ep

\bigskip


\section{Structure theorems for quasidirections, foldings and centralizers}

$\quad$ In this last section of the paper we will show that 
certain invariants (quasidirections, foldings, centralizers)
associated to elements of a given right-angled Artin group $(G, S)$ 
can be conveniently described in terms of the corresponding 
invariants associated to the so called {\em primitive} 
elements of $G$.

Before defining the primitive elements of $G$, we prove 
some useful statements on centralizers.

\begin{lem} Given $x, y, w \in G$ such that 
$x \bot y, x \subset w$ and $y \subset w$, the necessary 
and sufficient condition for $xy$ to belong to $Z_G(w)$ 
is that $x$ and $y$ belong to $Z_G (w)$.
\end{lem}

\bp The sufficiency part is trivial, so it 
remains to show that $x$ and $y$ belong to $Z_G (w)$ 
whenever $xy \in Z_G (w)$. By assumption 
$x \cup y = x \bullet y = y \bullet x$, and 
$w = x \bullet y \bullet z = z \bullet y \bullet x$, 
where $z = x^{-1} y^{-1} w$. We argue by induction on 
the length $d := l(w)$. Set $u := x \cap z, x^{\prime} := u^{-!} x$, 
and $z^{\prime} := u^{-1} z$. If $u = 1$ then  
$y \bullet z = z \bullet y$ by Lemma 3.9., and hence 
$x, y \in Z_G (w)$. Assuming that $u \neq 1$ and 
simplifying with $u$, we obtain 
$u^{-1} w = x^{\prime} \bullet y \bullet u \bullet  
z^{\prime} = z^{\prime} \bullet y \bullet u \bullet x^{\prime}$. 
As $x^{\prime} \bot z^{\prime}$ and $x^{\prime} \bot y$ 
thanks to \cite[Lemma 2.2.5.]{AC}, it follows by Lemma 3.9. again that 
$y \bullet x^{\prime} \bullet u = x^{\prime} \bullet y 
\bullet u = y \bullet u \bullet x^{\prime}$, therefore 
$x^{\prime} \bullet u = u \bullet x^{\prime}$. Since 
$y \bot u$, we get $u^{-1} w = x \bullet y \bullet z^{\prime} = 
z^{\prime} \bullet y \bullet x$, i.e. $xy \in Z_G(u^{-1} w)$. 
As $l(u^{-1}w) = d - l(u) < d$, it follows by the induction 
hypothesis that $y \in Z_G (u^{-1}w)$, and hence 
$x, y \in Z_G(w)$, since $yu = uy$. \ep

\begin{lem} For $x,y \in G$ and 
$m \geq  1, x^m = y^m$ implies $x = y$.
\end{lem}

\bp As $\varphi_x(1) = \varphi_{x^m} (1) = 
\varphi_{y^m} (1) = \varphi_y(1)$ by Lemma 6.8., we may 
assume without loss that $x$ and $y$ are both cyclically 
reduced. Setting $z :=x \cap y$ 
and assuming that $z \neq x$, let $s \in \widetilde{S}$ 
be such that $s \subset z^{-1} x$, 
whence $s \cap z^{-1} y = 1$. As 
$s \subset (z^{-1} x) \bullet x^{m - 1} = 
(z^{-1} y) \bullet y^{m - 1}$, it follows that 
$s \bot z^{-1}y$, and $s \subset y^{m - 1}$ by $(\perp)$. 
Since $\mathop{x \bullet x \bullet \ldots \bullet 
x}\limits_{\mathop{\underbrace{m \; {\mbox{factors}}}}} = 
\mathop{y \bullet y \bullet \ldots \bullet
y}\limits_{\mathop{\underbrace{m \; {\mbox{factors}}}}}$, 
the number of the $s\,'$s in any reduced 
decomposition of $x$ equals the corresponding number for 
$y$. As $s \subset z^{-1} x$ we necessarily have 
$s \sqsubset z^{-1}y$, contrary to $s \bot z^{-1} y$, 
by \cite[Lemma 2.2.5.]{AC} Consequently, $z = x \subset y$, 
and hence $x = y$ by symmetry. \ep

\begin{co} For $w \in G$ and 
$m \neq 0, Z_G (w^m) = Z_G(w)$.
\end{co}

\begin{pr} For any cyclically reduced 
element $w$ of $G$, the canonical isomorphism of median 
sets $\widetilde{1}^w \times X_{w, 1} 
\rightarrow X_w, (x,y) \mapsto x \bullet y = y \bullet x$ 
provided by {\em Corollary 7.5.} (with $a = 1 \in Z_G(w) \sse X_w$)
induces an isomorphism of median groups 
$\widetilde{1}^w \times H_w \rightarrow Z_G (w)$, 
where $H_w : = Z_G (w) \cap X_{w,1} = 
Z_G(w) \cap \mathop{1}\limits^{\equiv_w}$.
\end{pr}

\bp By Corollary 5.6., 
$  \mathop{1}\limits^{\sim_w} \subseteq Z_G (w)$ 
is a convex subgroup of $G$, so it is the {\em special}
subgroup of $G$ generated by $S_w := \{s \in S\,|\,s \perp w\}$.
It remains to show that 
$H_w$ is a median subgroup of $G$. First note that 
$H_w$ is a subgroup of $Z_G(w)$. Moreover the convex 
subset $\mathop{1}\limits^{\equiv_w}$ of $G$ is closed 
under the action from the left of $H_w$. Indeed, for 
$x \in Z_G(w)$ and $y \in G$, we get 
$y \mathop{\equiv}\limits_w 1 \Llra 
x^{-1} y \mathop{\equiv}\limits_{x^{-1} wx} x^{-1} \Llra
x^{-1} y \mathop{\equiv}\limits_w x^{-1}$, and 
$x^{-1} \mathop{\equiv}\limits_w 1 \Llra 
1 \mathop{\equiv}\limits_{x wx^{-1}} x \Llra 
1 \mathop{\equiv}\limits_w  x$, as required.

Thus it remains to show that $Y(x, y, z) \in Z_G(w)$ 
for all $x, y, z \in H_w$. Since 
$X_{w, 1} = \mathop{\cup}\limits_{n \geq 0} [w^{-n} , w^n]$
by Proposition 7.3.$(1)$, it follows that
there is $n \geq 1$ such that $x, y, z \in 
[w^{-n}, w^n] \subseteq X_{w, 1}$. 
Consequently, the elements $w^n x, w^ny$ and $w^nz$ belong 
to $Z_G(w) \cap [1, w^{2n}] \subseteq Z_G(w^{2n}) \cap [1, w^{2n}]$, 
therefore $w^n Y(x, y, z) = Y(w^n x, w^n y, w^n z) \in Z_G(w^{2n})$ 
by Lemma 3.10. As $Z_G(w^{2n}) = Z_G(w)$ according to Corollary 8.3., 
we conclude that $Y(x, y, z) \in Z_G(w)$ as desired. \ep

\begin{lem} Let $w$ be a cyclically reduced element 
of $G$, and $a \in  X_w$, so the conjugate $a^{- 1} w a$ of
$w$ is cyclically reduced too. Then, the inner group automorphism 
$x \mapsto a^{-1} x a$ of the group $G$ induces an isomorphism of 
median groups $H_w \rightarrow H_{a^{-1} w a}$.
\end{lem}

\bp By Corollary 7.5., the convex subsets 
$\widetilde{1}^w$ and $X_{w, 1}$ are orthogonal, i.e. 
$x \bot y$ provided $x \in \widetilde{1}^w$ and 
$y \in X_{w, 1}$, therefore the mapping $x \mapsto a^{-1} xa$ 
is the identity on $H_w$ whenever $a \in \widetilde{1}^w$. 
Thus we may assume without loss that $a\in X_{w, 1}$. 
First notice that $a^{-1} H_w a = H_{a^{-1} wa} = 
Z_G(a^{-1}w a) \cap a^{-1}X_{w, 1}$. Indeed, for $x \in 
H_w$ we get $\varphi_w(xa) = x \varphi_w (a) = xa$, i.e. 
$xa \in X_w$. On the other hand, 
$a \mathop{\equiv}\limits_{w} 1$ implies 
$xa \mathop{\equiv}\limits_{xwx^{-1}} x$, therefore 
$xa \mathop{\equiv}\limits_{w} 1$, so $xa \in X_{w, 1}$, 
since $wx = xw$ and $x \mathop{\equiv}\limits_{w} 1$. 

As $a \in X_{w, 1}$ and $b \subset a \Lra
b^{-1}a \in b^{-1} X_{w, 1} = b^{-1} X_{w, b} = 
X_{b^{-1}  w b, 1}$, proceeding by induction on the length 
$d := l(a)$, we are reduced to the case $d = 1$, i.e. 
$a \in \widetilde{S} \cap X_{w, 1}$. Since for each
$x \in H_w \subseteq X_{w, 1}$ there is $n \geq 0$ such 
that $[1, x]=[w^{-n} \cap x, w^n \cap x]$, and 
$w^m \cap x \in H_w$ for all $m \in \Z$ by Proposition 
8.4., and since $H_{w^m} = H_w$ for all $m \neq 0$ by 
Corollary 8.3, it remains to show that 
$a^{-1} xa \subset a^{-1} ya$ whenever $x \subset y 
\subset w$. As $a \in \widetilde{S} \cap X_{w, 1}$ 
we distinguish the following two cases : 

Case $(1)$ : $a \subset w$. If $a \cap x = 1$, whence
$a^{-1} \bot x$, then we are done by Lemma 4.10., so we may
assume that $a \subset x$. Setting 
$x^{\prime } := a^{-1} x, y^{\prime } := x^{-1} y$ and 
$z := y^{-1} w$, we get 
$w = a \bullet x^{\prime } \bullet y^{\prime } \bullet z = 
z \bullet a \bullet x^{\prime} \bullet y^{\prime } = 
y^{\prime } \bullet z \bullet a \bullet x^{\prime}$, therefore 
$a \subset z \bullet a$ by $(A_1)$, and hence 
$a \subset y^{\prime } \bullet a$ by $(A_1)$ again. As 
$w^2 = w \bullet w$, it follows that 
$a^{-1} x a = x^{\prime } \bullet a \subset x^{\prime } 
\bullet y^{\prime } \bullet a = a^{-1} y a$ as required.

Case $(2)$ : $a \subset w^{-1}$. As 
$w^{-1} y \subset w^{-1} x \subset w^{-1}$, we may apply 
Case $(1)$ to get 
$a^{-1} w^{-1} y a  \subset a^{-1}w^{-1}xa \subset a^{-1} w^{-1} a$, 
and hence $a^{-1} x a \subset a^{-1} y a$, as desired. \ep

\begin{rem} \em 
For a cyclically reduced element $w$ of $G$ 
and an element $a \in X_w \backslash X_{w, 1}$, the group isomorphism 
$Z_G(w) \rightarrow Z_G(a^{-1} w a), x \mapsto a^{-1} x a$, 
is not necessarily an isomorphism of median groups. For 
instance, let $S = \{a,b,c\}$, and let $G$ be 
given by the presentation 
$G = < S; [a,c] = [b, c] = 1>$, so 
$G \cong F_2 \times \Z$. We obtain 
$Z_G (c) = X_c = G, X_{c,1} = H_c = < c>$, and $a \subset ab$, 
but $a^{-1} aa=a \not\subset ba = a^{-1} (ab)a$.
\end{rem}

The following definition is justified by Lemma 8.5.

\begin{de} A non-trivial element $w$ of $G$ 
is called {\em primitive} if for some (for all) $a \in X_w$, 
the median subgroup 
$H_{a^{-1}wa} = Z_G(a^{-1} w a) \cap X_{a^{-1} wa, 1}$ 
is cyclic, generated by $a^{-1} wa$. 
\end{de}

In particular, a cyclically 
reduced element $w \neq 1$ is primitive iff $H_w$ is 
generated by $w$.
As $X_{xwx^{-1}} = x X_w$ for all $x, w \in G$, 
the primitiveness is preserved by conjugation.

The next lemma provides equivalent descriptions for 
primitive elements.

\begin{lem} The following assertions are 
equivalent for an element $w \neq 1$ of $G$.

$(1)\, w$ is primitive.
 
$(2)\,$ The cell $C := [1, \varphi_w (1)^{-1} w \varphi_w (1)]$ 
is quasilinear, i.e. $|\partial C| \leq 2$, and $w$ is not a 
proper power of some element of $G$.

$(3)\,$ For all $x \in G$, the cell $[1, xwx^{-1}]$ is 
quasilinear, and $w$ is not a proper power of some 
element of $G$.
\end{lem}

\bp $(1) \Lra (3)\,$. Assume that $w$ is 
primitive, and let $x \in G$. By Lemma 4.2., 
$xwx^{-1} = u \bullet v \bullet u^{-1}$, where 
$u = \varphi_{xwx^{-1}} (1) = x \varphi_w (x^{-1})$, 
and $v = u^{-1}xwx^{-1} u$ is cyclically reduced and primitive. 
In particular, $v$ (and hence $w$) cannot be a proper power 
since assuming $v = v^{\prime n}$ for some $v^{\prime} \in G, n \geq 1$, 
we get $v^{\prime} \in H_{v^{\prime}} = H_{v}$, and hence 
$n = \pm 1$. Assuming that $[1, x wx^{-1}] = [a, b]$ for some 
$a, b \in G$, i.e. $a \bot b$ and 
$x w x^{-1} = a \bullet b = b \bullet a$, and setting 
$u_1 := u \cap a, u_2 := u \cap b, a^{\prime} := u^{-1}_{1} a u_1$ 
and $b^{\prime } := u_{2}^{-1} b u_2$, we obtain 
$a = u_1 \bullet a^{\prime } \bullet u_{1}^{-1}, b = 
u_2 \bullet b^{\prime } \bullet u_{2}^{-1}$, and 
$[1,v] = [a^{\prime }, b^{\prime }]$, therefore either 
$a^{\prime } = 1 $ or $b^{\prime } = 1$ since 
$a^{\prime }, b^{\prime} \in [1, v] \cap Z_G (v) = \{1,v\}$. 
In the former case we get $a = 1$, while in the latter 
case it follows that $b = 1$. Thus the cell 
$[1, x w x^{-1}]$ is quasilinear as required. 

$(3) \Lra (2)\,$ is trivial.

$(2) \Lra (1)\,$ It suffices to show that 
$[1, w] \cap Z_G (w)=\{1, w\}$ whenever the cyclically 
reduced element $w$ is not a proper power and the cell 
$[1, w]$ is quasilinear. Let $a \in [1,w] \cap Z_G (w)$ 
be such that $a \neq 1$ and its length $l(a)$ is minimal. 
We have to show that $a = w$. Let $n \geq 1$ be such that 
$a^n = \mathop{\underbrace{a \bullet a \bullet \ldots \bullet
a}}\limits_{n \; \; {\mbox{factors}}} \subset w$ and 
$a^{n + 1} \not\subset w$. Setting $b := a^{-n} w$, we obtain 
$w = a \bullet b = b \bullet a $. As $Z_G(w)$ is a median 
subgroup of $G$ we get $a \cap b \in [1, w] \cap Z_G (w)$, 
therefore $a \cap b =1$ by the minimality of $l(a)$ and the 
maximality of $n$. Consequently, $[1, w] = [a^n, b]$, and hence 
$w = a^n$ since the cell $[1, w]$ is quasilinear by 
assumption. As $w$ is not a proper power, we obtain 
$a = w$ as desired. \ep

Let $Prim (G)$ denote the subset of all primitive 
elements  of $G$. Obviously, 
$\widetilde{S} \subseteq  Prim (G)$, and $Prim (G)$ 
is closed under the operation $w \mapsto w^{-1}$. In 
particular, if $G$ is freely generated by $S$ then 
$Prim (G)$ consists of those $w \in G \backslash \{1\}$ 
which are not proper powers, while 
$Prim (G) =  \widetilde{S}$ whether $G$ is the free Abelian 
group generated by $S$.  

The elements of $G$ admit canonical representations as 
products of powers of commuting primitive elements, as 
follows.

\begin{te} For a given element $w \in G$, there exist 
primitive elements $p_1, \ldots , p_n$ and positive integers 
$m_1, \ldots , m_n$ such that $a^{-1} p_i a \bot a^{-1} p_ja $ 
for $i \neq j$ and $a \in X_w$ (in particular, the 
$p_{i}\,^{\prime }s$ are commuting primitive elements), and 
$w = \prod\limits_{i = 1}^{n} p_{i}^{m_i}$. The pairs 
$(p_i, m_i)$ are uniquely determined up to a permutation 
of the indices $i = 1, \ldots , n$.
\end{te} 

\bp For all 
$a \in X_w, \; [1, a^{-1} wa] \cap Z_G (a^{-1} wa)$ 
is a median subset of the median group 
$H_{a^{-1} wa} = Z_G(a^{-1} wa) \cap a^{-1} X_{w, a}$. Given 
$a \in X_w$, let $u_1, \ldots ,u_n$ be the minimal elements of 
$[1, a^{-1} wa] \cap Z_G \; (a^{-1} wa)$ with respect 
to the order $\subset $. Obviously, the 
$u_{i}^{\prime }s$ are cyclically reduced and pairwise 
orthogonal. Also they are not proper powers by Corollary 8.3. 
Moreover, assuming $[1, u_i]= [u_{i}^{\prime }, 
u_{i}^{\prime \prime}]$, it follows by Lemma 8.1. that 
$u_{i}^{\prime } \in [1,a^{-1} wa] \cap Z_G (a^{-1} wa)$, 
therefore $u_{i}^{\prime } \in \{1, u_i\}$ by the minimality of 
$u_i$. Thus the cells $[1, u_i], i = 1, \ldots , n$, are 
quasilinear, and hence the $u_{i}^{\prime }s$ are 
primitive according to Lemma 8.8. Let $m_i \geq 1$ be 
the largest natural number for which 
$u_{i}^{m_i} \subset a^{-1} wa$, and let $u = 
\displaystyle\cup_{i = 1}^{n} u_{i}^{m_i} = 
\prod\limits_{i = 1}^{n} u_{i}^{m_i}$ and 
$v = u^{-1} (a^{-1} wa)$. As 
$u \in [1, a^{-1} w a] \cap Z_G (a^{-1} wa)$, we get 
$v \in [1, a^{-1} wa] \cap Z_G (a^{-1} wa)$, therefore, 
assuming $v \neq 1$, there is an index $i$ such that 
$u_i \subset w$. Writing $u = u_{i}^{m_i} \bullet u^{\prime }$, 
with $u_i \bot u^{\prime }$, and $v = u_i \bullet v^{\prime }$, 
we get $a^{-1} wa = u \bullet v=u_{i}^{m_i} \bullet u^{\prime} 
\bullet u_i \bullet v^{\prime }= 
u_{i}^{m_i + 1}  \bullet u^{\prime } \bullet v^{\prime }$, 
contrary to the definition of $m_i$. Consequently, 
$a^{-1} wa = u$. Setting $p_i := a u_i a^{-1} $ for 
$i = 1, \ldots n$, we obtain a representation of 
$w = \prod\limits_{i = 1}^{n} p_{i}^{m_i}$ as a 
product of powers of the commuting primitive elements 
$p_1, \ldots , p_n$, so to end the proof of the existence 
part of the statement, it remains to note that 
$b^{-1} p_i b \bot b^{-1} p_j b$ for $i \neq j$ and 
$b \in X_w$, since the conjugation map 
$x \mapsto (a^{-1} b)^{-1} x(a^{-1}b)$, where 
$a^{-1} b \in a^{-1} X_w = X_{a^{-1} w a}$, induces an 
isomorphism of median groups 
$H_{a^{-1} wa} \longrightarrow H_{b^{-1} wb}$ 
according to Lemma 8.5.

To prove the uniqueness up to permutation, assume that 
the pairs $(p_{i}, m_i), i = 1, \ldots, n$, satisfy the requirements of the 
statement. If suffices to show that for all $a \in X_w$, 
the $a^{-1} p_i a\,^{\prime }s $ are minimal elements of 
the lattice $L: = [1, a^{-1} wa] \cap Z_G(a^{-1} wa)$. As 
$u_i: = a^{-1} p_i a \neq 1 $ belongs to the lattice $L$ by 
assumption, there is a minimal element $v$ of $L$ such that 
$v \subset u_i$. Since $u_i^{m_i} \bot u_i^{-m_i} (a^{-1} wa)$ 
and $v \in Z_G (a^{-1} wa)$, we get $v \in Z_G( u_i^{m_i})$, 
therefore $v \in [1, u_i]\cap Z_G (u_i)$ by Corollary 8.3., 
and hence $v = u_i$ since $v \neq 1$ and $u_i$ is cyclically 
reduced and primitive by assumption. \ep

\begin{rem} \em 
A similar result with Theorem 8.9. above
is proved in \cite{Bau1, Serv} by different methods.
\end{rem}

For any $w \in G$, let $Prim (w) \subseteq  Prim (G)$ 
denote the finite set $\{p_1, \ldots , p_n\}$ of primitive 
elements uniquely associated to $w$ by Theorem 8.9. Let 
$Prim (w)^{\sim }$ denote the disjoint union of 
$Prim (w)$ and $Prim (w^{-1}) = Prim (w)^{-1}$. 
Notice that $Prim (xwx^{-1}) = x  Prim  (w) 
x^{-1}$ and $Prim (xwx^{-1})^{\sim} = x  
Prim  (w)^{\sim } x^{-1}$ for all $x \in G$.

For any $w \in G$, we have denoted by $S_w$ the subset of $S$ 
consisting of those $s \in S$ for which $s \bot w$.
In particular, $S_t = \{s \in S \sm \{t\}\,|\,st = ts\}$
for all $t \in S$.

The next statement is an immediate consequence of 
Proposition 8.4. and Theorem 8.9.

\begin{co} $(1)\,$ For any cyclically reduced 
element $w$ of $G, S_w = \bigcap\limits_{p \in Prim  (w)} S_p$ 
generates the convex subgroup $\widetilde 1^w
= \mathop{\bigcap}\limits_{p \in \; Prim  (w)} \widetilde 1^p$, 
while the median subgroup $H_w = Z_G (w) \cap X_{w, 1}$ is Abelian, 
freely generated by $Prim  (w)$ and contained in the center of 
$Z_G(w)$.

$(2)\,$ For any $w \in G$ and $a \in X_w, Z_G(w)$ is the direct 
product of the right-angled Artin group generated by $\bigcap\limits_{p \in 
Prim  (w)}  a S_{a^{-1} pa} a^{-1}$ and the free Abelian group 
generated by $Prim  (w)$. In particular, $Z_G (w)$ is a right-angled
Artin group.

$(3)\,$ The center $Z(G)$ of $G$ is the free Abelian group 
generated by the (possibly empty) set 
$\{s \in S\,|\,\forall t \in S, st = ts\}$.
\end{co}

\begin{rems} \em
$(1)\,$ It is known \cite{Bau2} that, by contrast with 
free groups and free Abelian groups, in general,
the partially commutative freeness is not transferable to 
arbitrary subgroups. A graph theoretic transfer criterion
for right-angled Artin groups $(G, S),\,S$ finite, is 
given in \cite{Dr}.  

$(2)\,$ A result similar with Corollary 8.11. is proved
in \cite{Bau1, Serv} by different methods.
\end{rems}

\begin{co} Let $w \in G$. Then, the following assertions hold.

$(1)\,$ For all 
$a \in X_w, X_{w, a}$ is the convex closure of the union of 
its convex subsets $X_{p, a}$ for $ p \in Prim (w)$, and 
$X_{w, a} \cong \prod\limits_{p \in  Prim  (w)} X_{p, a}$.

$(2)\,$ $X_w =  \displaystyle\cap_{p \in  Prim (w)} X_p$.

$(3)\,$ The folding $\varphi_w$ is obtained by composing 
the commuting foldings $\varphi_p$ for $p \in  Prim (w)$.

$(4)\,$ For $u \in G, \varphi_u \leq \varphi_w$, i.e. 
$X_u \subseteq X_w$, provided
$Prim  (w)^{\sim} \subseteq  Prim (u)^{\sim}$.
\end{co}

\bp $(1)\,$ The inclusion $X_w \subseteq X_p$ for 
$p \in Prim (w)$ is immediate by Theorem 8.9. Writing 
$w = \prod\limits_{p \in \; Prim (w)} p ^{m_p}$ with 
$m_p \geq 1$, we get easily 
$$[w^{-n}a, w^n a] =
[\displaystyle\cup_{p \in \; Prim  (w)} [p^{-nm_p} a, 
p^{nm_p} a]] \cong \prod\limits_{p \in \; Prim  (w)}[p^{-nm_p} 
a, p^{nm_p} a],$$ 
therefore
$X_{w,a} = [\displaystyle\cup_{p \in  Prim  (w)} X_{p,a}] 
\cong \prod\limits_{p \in  Prim  (w)} X_{p, a}\,$ 
for all $a \in X_w$.

$(2)\,$ Given $x \in \displaystyle\cap_{p \in  Prim  (w)} X_p$ 
and $a \in X_w$, it follows by Corollary 7.5. that 
$[a, x] = [y_p, z_p]$ with $y_p \in \widetilde a^p$ and 
$z_p \in X_{p, a}$ for $p\in  Prim  (w)$, and hence 
$[a, x] = [y, z]$, where 
$y = \mathop{\vee}\limits_a \{y_p\,|\, p \in Prim  (w)\}$ and 
$z = \mathop{\vee}\limits_x \{z_p\,| \; p \in Prim (w) \}$, 
since the negation operator $\neg$ is a median set 
automorphism of $\partial [a, x]$. As 
$y \in \mathop{\cap}\limits_{p \in Prim  (w) }
\widetilde a^p =  \widetilde{a}^w$ and $z \in 
[\mathop{\cup}\limits_{p \in 
 Prim  (w)} X_{p, a}] = X_{w, a}$, we obtain 
$x \in [y, z] \subseteq [\widetilde{a}^w \cup X_{w, a}] = X_w $
as desired.

$(3)\,$ and $(4)\,$ are immediate  consequences of $(2)$. \ep

\begin{rem} \em 
The converse of the assertion $(4)\,$ above is 
not necessarily true. For instance, if 
$G = <s, t; [s, t] = 1> \cong \Z \times \Z$ then 
$\varphi_s = \varphi_t = 1_G$, but 
$Prim  (s)^{\sim } = \{s, s^{-1} \} 
\neq \{t, t^{-1}\} = Prim  (t)^{\sim }$.
\end{rem}

\begin{co} For all $w \in G, 
\mathop{\bullet}\limits_{w} = 
\mathop{\cap}\limits_{p \in Prim  (w)} 
\mathop{\bullet}\limits_{p} $, i.e. 
for all $x, y \in G, x \mathop{\preceq }\limits_{w } y \Llra 
x \mathop{\preceq }\limits_{p } y$ for all $p \in Prim  (w)$.
\end{co}

\bp Proceeding by induction on the distance 
$d := d(x, y)$ we are reduced to the case $d = 1$, i.e. 
$x^{-1} y = s \in \widetilde{S}$. Without loss we may 
also assume that $x = 1$ and $y = s \in \widetilde{S}$. 
Thus we have to show that $s \subset ws \Llra \forall 
p \in Prim (w),\,s \subset ps$. Setting $a := \varphi_w (1)$, 
we get $w = a \bullet w^{\prime} \bullet a^{-1}$, where 
$w^{\prime } := a^{-1} wa$ is cyclically reduced. 
By Lemma 6.8. and Theorem 8.9. it suffices to show that 
for $u^{\prime }, v^{\prime } \subset w^{\prime } $ such that 
$[1, w^{\prime }] = [u^{\prime }, v^{\prime }]\,$, $s \subset ws 
\Llra s \subset us$ and $s \subset vs$, where 
$u = a u^{\prime } a^{-1}$ and $v = a v^{\prime } a^{-1}$. As 
$\varphi_{w^{\prime }} = \varphi_{u^{\prime }} \circ 
\varphi_{v^{\prime }} = \varphi_{v^{\prime}} \circ 
\varphi_{u^{\prime }}$ by Corollary 8.13.$(3)$, we get 
$\varphi_{w^{\prime}} (a^{-1}) = 1 \in [\varphi_{u^{\prime }} 
(a^{-1}), \varphi_{v^{\prime }} (a^{-1})] \subseteq [1, a^{-1}]$, 
and hence $\varphi_{u^{\prime }} (a^{-1}) \bot 
\varphi_{v^{\prime }} (a^{-1})$.

Setting $b := \varphi_{u^{\prime }} (a^{-1})^{-1}, \; c := 
\varphi_{v^{\prime }} (a^{-1})^{-1}$ and 
$a^{\prime } := ab^{-1} c^{-1}$, it follows that $b 
\bot c, \; a = a^{\prime } \bullet b \bullet c = 
a^{\prime } \bullet c \bullet b, \; \varphi_u (1) = a 
\varphi_{u^{\prime }} (a^{-1}) = a^{\prime } \bullet c$, and 
$\varphi_v (1) = a \varphi_{v^{\prime }} (a^{-1}) = a^{\prime } 
\bullet b$. As $u \cap a^{-1} = u^{-1} \cap a^{-1} = 1$, 
we also obtain $\varphi_{u^{\prime }} (a^{-1}) = u^{\prime } 
\bullet a^{-1} \cap a^{-1}= u^{\prime -1} \bullet 
a^{-1} \cap a^{-1} = u^{\prime } \bullet a^{-1} 
\cap u^{\prime -1} \bullet a^{-1}$, in particular, 
$b \bot u^{\prime }$, and similarly $c \bot v^{\prime }$. 
Setting $w^{\prime \prime } := a^{\prime -1} w 
a^{\prime }, \; u^{\prime \prime } := 
a^{\prime -1} u a^{\prime }$ and $v^{\prime \prime } :=
a^{\prime -1} v a^{\prime }$, it follows that $[1, 
w^{\prime \prime }] = [u^{\prime \prime }, v^{\prime \prime }], w= 
a^{\prime } \bullet w^{\prime \prime }  \bullet a^{\prime -1} ,u= 
a^{\prime } \bullet u^{\prime \prime } \bullet a^{\prime -1}$ and $v= 
a^{\prime } \bullet v^{\prime \prime } \bullet a^{\prime -1}$. We 
distinguish the following two cases :

Case $(1)\,$ : $ws = w \bullet s = a^{\prime } \bullet w^{
\prime \prime } \bullet a^{\prime -1} \bullet s$. Then, 
$u^{-1} \cap s = a^{\prime } 
\bullet u^{\prime \prime -1} \bullet a^{\prime -1} 
\cap s = a^{\prime } 
\bullet u^{\prime \prime -1} \cap s \subset w^{-1} \cap s = 1$, 
and, similarly, $v^{-1} \cap s = 1$, so $us = u \bullet s$ 
and $vs = v \bullet s$. Assuming that $s \subset ws 
= a^{\prime } \bullet  u^{\prime \prime }  
\bullet  v^{\prime \prime }  
\bullet  a^{\prime -1}  \bullet s$, but 
$s \not\subset us =  a^{\prime } 
\bullet  u^{\prime \prime }  \bullet  a^{\prime -1}  \bullet 
s$, we get $s \bot a^{\prime } \bullet u^{\prime \prime }$ 
and hence $us = s \bullet a^{\prime } 
\bullet  u^{\prime \prime } \bullet a^{\prime -1}$, 
contrary to the assumption $s \not\subset us$. 
Consequently, $s \subset us$ and $s \subset vs$ whenever 
$s \subset ws$. Conversely, assuming $s \subset us$,
but $s \not\subset ws$, we get 
$s \bot a^{\prime } \bullet u^{\prime \prime }$, whence 
$a^{\prime } \bullet  v^{\prime \prime } \bullet s \subset ws$. 
It follows that $s \cap a^{\prime } \bullet  
v^{\prime \prime } \bullet s = 1$, and hence 
$s \not\subset vs$, as required.

Case $(2)\,$ : $s \subset w^{-1} = a^{\prime } \bullet  w^{\prime 
\prime -1} \bullet a^{\prime -1}$, whence $s \subset a^{\prime } 
\bullet  w^{\prime \prime -1}$. If $s \subset a^{\prime }$ then we 
have nothing to prove, so let us assume that $s \bot a^{\prime }$ and 
$s \subset  w^{\prime \prime -1}$. As $[1,  w^{\prime \prime -1}]=[ 
u^{\prime \prime -1},  v^{\prime \prime -1}]$, we may assume that $s 
\subset  u^{\prime \prime -1}$ and $s \bot v^{\prime \prime }$. 
Thus $s \bot v$, therefore $s \subset s \bullet v = v \bullet s$. 
Setting $ u^{\prime \prime } = u^{\prime \prime \prime} 
\bullet s^{-1}$, we get $us = 
a^{\prime } \bullet  u^{\prime \prime \prime } 
\bullet a^{\prime -1}$ and $ws = a^{\prime } 
\bullet u^{\prime \prime \prime }\bullet 
v^{\prime \prime}\bullet a^{\prime -1}$, therefore 
$s \subset us \Llra s \subset   u^{\prime \prime \prime } 
\Llra s \subset ws $ as desired. \ep

The next statement provides a classification of the 
quasidirections $\mathop{\bullet}\limits_{w}$ for $w \in G$.

\begin{pr} The mapping 
$w \in G \mapsto Prim (w)$ induces an antiisomorphism 
of the ordered set of the quasidirections 
$\mathop{\bullet}\limits_{w}$ for $w \in G$ onto the set 
${\cal F} (G)$ ordered by inclusion, consisting of 
the finite subsets $F \subseteq \; Prim (G)$ satisfying

$(i)\,X_p \cap X_q \neq \emptyset $ for $ p, q \in F$, and

$(ii)\,a^{-1} pa \bot a^{-1} qa$ for $p, q \in F, p \neq q$, 
and for some (for all) $a \in X_p \cap X_q$.
\end{pr}

\bp Notice that for any finite subset 
$F \subseteq \; Prim \; (G)$, the conditions $(i)$ and $(ii)$ 
above are equivalent with the apparently stronger 
conditions

$(i^{\prime }) \mathop{\cap}\limits_{p \in F}  X_p \neq \emptyset $, 
and 

$(ii^{\prime }) \; a^{-1} pa \bot a^{-1}qa$ for 
$p,q \in F, p \neq q$, and for some (for all) 
$a \in \mathop{\cap}\limits_{u \in F} X_u$.

According to Corollary 8.15., it remains to show that 
for $w \in G$ and $p \in  Prim  (G)$, $p \in  Prim  (w)$ 
whenever the preorder $\mathop{\preceq}\limits_p$ is finer 
than the preorder  $\mathop{\preceq}\limits_w$. Without 
loss we may assume that $w$ is cyclically reduced, i.e. $1 \in X_w$. 
First let us show that $X_w \cap X_p \neq \emptyset $. 
Since $\varphi_p (1) \mathop{\preceq}\limits_w 
\varphi_w (\varphi_p (1))$ by Lemma 6.9., it follows by 
assumption that 
$\varphi_p (1) \mathop{\preceq}\limits_p \varphi_w (\varphi_p (1))$. 
On the other hand, as 
$\varphi_w (\varphi_p (1)) \in [1, \varphi_p (1)]$ and 
$1 \mathop{\ll}\limits_p \varphi_p (1)$ by Lemma 6.9. again, we get 
$\varphi_p(1) = \varphi_w (\varphi_p (1)) \in X_w \cap X_p$
as required. Next let us show that $wx \in X_p$ whenever 
$x \in X_w \cap X_p$. As $x \in X_p \Lra
\varphi_p (wx) \in [x, wx]$, and $x \in X_w \Lra 
x \mathop{\ll}\limits_w wx$, it follows that 
$\varphi_p (wx) \mathop{\ll}\limits_w wx$, whence 
$\varphi_p (wx) \mathop{\preceq}\limits_p wx$. Since 
$wx \mathop{\ll }\limits_p \varphi_p (wx)$ by Lemma 6.9., 
we get $wx = \varphi_p (wx) \in X_p$ as desired. Consequently, 
$w^n x \in X_p$ for all $n \geq 0$ provided $x \in X_w \cap X_p$. 
Thus for all $x \in X_w \cap X_p$ and for all $n \geq 0$, 
the element $(w^n x)^{-1} p (w^n x)$ is a cyclically 
reduced conjugate of $p$. 
Since there are only finitely many cyclically reduced 
conjugates of $p$, there is some $n \geq 1$ such that 
$p \in Z_G (w^n)$, and hence $p \in Z_G (w)$ by Corollary 8.3. 
As $p$ is primitive, the cell $[1, p]$ is quasilinear 
by Lemma 8.8., and hence, according to Proposition 8.4., 
either $p \in \widetilde{1}^w $ or  
$p \in H_w = Z_G (w) \cap X_{w, 1}$. The former case 
would imply $p \in \widetilde{1}^p$, i.e. $p = 1$, a 
contradiction, so $p \in H_w$, whence 
$p \in \; Prim (w)^{\sim}$. The assumption 
$p^{-1} \in \; Prim \; (w)$ would imply 
$1 \mathop{\ll}\limits_{w} p^{-1}$, whence 
$1 \mathop{\preceq }\limits_p p^{-1}$, i.e. $p = 1$, 
again a contradiction. Consequently, $p \in \; Prim \; (w)$ 
as required. \ep

\begin{co} The mapping 
$(w,a) \in G \times G \mapsto 
(\stackrel{\equiv}{a^p})_{p \in \; Prim \; (w)}$ induces 
a bijection of the set of directions 
$\{\mathop{\vee}\limits_{w; a}\,|\,(w, a) \in G \times G\}$ 
onto the disjoint union 
$\mathop{\bigsqcup}\limits_{F \in {\cal F} (G)} 
G/\mathop{\bigcap}\limits_{p \in F} \equiv _p$.
\end{co}

\bp By Propositions 7.9. and 8.16., it suffices 
to show that for 
all $w, u, a, b \in G$, the quasidirections 
$\mathop{\bullet}\limits_{w}$ and $\mathop{\bullet}\limits_{u}$
coincide whenever the directions $\mathop{\vee}\limits_{w; a}$ 
and $\mathop{\vee}\limits_{u; b}$ coincide, since then we get 
$\mathop{\vee}\limits_{w; a} = \mathop{\vee}\limits_{u,b} \Llra 
\mathop{\bullet}\limits_{w} = \mathop{\bullet}\limits_{u}$ 
and $a \mathop{\equiv}\limits_w b \Llra 
Prim (w) = Prim (u)$ and $a \mathop{\equiv}\limits_p b$ 
for all $p \in Prim (w)$. Let $w, u, a, b \in G$ be such 
that $\mathop{\vee}\limits_{w; a} = 
\mathop{\vee}\limits_{u; b}$. As 
$\mathop{\vee}\limits_{w; a} = \mathop{\vee}\limits_{w; \varphi_w (a)}$ 
we may assume from the beginning that $a \in X_w$ and 
$b \in X_u$. Moreover we may assume that 
$a = b \in X_w \cap X_u$ since $a  
\mathop{\ll}\limits_w a \mathop{\bullet}\limits_w b = a 
\mathop{\vee}\limits_{w; a} b = a \mathop{\vee}\limits_{u; a} b = b 
\mathop{\bullet}\limits_u a \mathop{\gg}\limits_u b$. 
Of course we may also assume that $a=b=1 \in X_w \cap X_u$, 
so $\{x \in X_{w, 1} \mid 1 \mathop{\ll}\limits_w x \} = 
\{x \in X_{u,1}   \mid 1 \mathop{\ll}\limits_u x \}$. 
In particular, $1 \mathop{\ll}\limits_w u^n \in X_{w, 1}$ 
for all $n \geq 0$. Since there are only finitely many 
cyclically reduced conjugates of $w$ it follows that $u$ 
belongs to the cone of positive elements with respect to 
the order $\mathop{\ll}\limits_w$ of the free abelian 
group $H_w = Z_G (w) \cap X_{w, 1}$ generated by $Prim \; (w)$. 
Consequently, $Prim \; (u) \subseteq \; Prim \; (w)$, 
and hence, by symmetry, $Prim \; (w) = Prim \; (u)$, 
therefore $\mathop{\bullet}\limits_{w} 
= \mathop{\bullet}\limits_{u}$ by Proposition 8.16. \ep

\begin{co} For all $w \in G, \; Stab \; 
(\mathop{\bullet}\limits_{w}) = \; Stab \; (\varphi_w) = 
Stab \; (X_w) = Z_G 
(w) = \mathop{\bigcap}\limits_{p \in \; Prim \; (w)} Z_G (p)$, 
and for 
all $a \in G, \; Stab \; ( \mathop{\vee}\limits_{w; a}) = Stab \; 
(\Psi_{w, \varphi_w (a)} ) = Stab \; (X_{w, \varphi_w (a)})$ 
is the free Abelian group generated by $Prim  (w)$.
\end{co}

\bp The inclusions 
$\mathop{\bigcap}\limits_{p \in Prim \; (w)} 
Z_G (p) \subseteq Z_G (w) \subseteq \; Stab ( 
\mathop{\bullet}\limits_{w}) \cap 
\; Stab \; (\varphi_w)$ are trivial. Assuming $x \in Stab \; 
(\mathop{\bullet}\limits_{w})$, i.e. 
$\mathop{\bullet}\limits_{xwx^{-1}} = 
\mathop{\bullet}\limits_{w}$, it follows by Proposition 8.16. 
that $Prim  (xwx^{-1}) = x \; Prim  (w) x^{-1} = Prim 
 (w)$, so $x^n p x^{-n} = p$ for $p \in \; Prim  (w)$ 
and a divisor $n$ of $\mid Prim (w)\mid !$, therefore 
$ p \in Z_G (x^n) = Z_G (x)$ for all $ p \in \; Prim (w)$, i.e. 
$x \in \mathop{\bigcap}\limits_{p \in \; Prim (w)} \; Z_G 
(p)$ as desired. Assuming that
$x \in \; Stab \; (\varphi_w) = \; Stab \; 
(X_w)$, i.e. $x X_w = X_w$, and taking some 
$a \in X_w$, we obtain a family 
$(a^{-1} x^{-n} w x^n a)_{n \in \Z}$ of cyclically 
reduced conjugates of $w$, and hence $x \in Z_G (w)$ 
by the finiteness argument used in the proof of Corollary 8.17..

On the other hand, it follows by Corollary 8.17. that 
$Stab (\mathop{\vee}\limits_{w; a}) = 
\{ x \in G \mid \mathop{\vee}\limits_{xwx^{-1}; xa} = 
\mathop{\vee}\limits_{w; a}\} = \{ 
x \in Z_G (w) \mid \; x \varphi_w (a) 
\mathop{\equiv}\limits_{w} \varphi_w 
(a)\}=$ the free abelian group generated by 
$Prim \; (w)$. We get a similar result for 
$Stab \; (\Psi_{w, \varphi_w (a)})$ since one 
checks easily, as in the proof of Corollary 8.17., 
that for all $w, u \in G, a \in X_w$ and 
$b \in X_u, \Psi_{w, a} = \Psi_{u, b} \Llra 
Prim  (w) ^{\sim} = Prim  (u)^{\sim}$ and 
$a \mathop{\equiv}\limits_{w} b$. \ep

\end{document}